%% file: main.tex
\documentclass[review,onefignum,onetabnum]{siamart250211}
\input{ex_shared}

\ifpdf
\hypersetup{
  pdftitle={Parameter-robust preconditioners for a cell-by-cell poroelasticity model with interface coupling},
  pdfauthor={M. Causemann and M. Kuchta}
}
\fi

\begin{document}

\nolinenumbers
\maketitle

% REQUIRED
\begin{abstract}
  This paper presents a scalable and robust solver for a cell-by- cell poroelasticity model, describing the mechanical interactions between brain cells embedded in extracellular space. Explicitly representing the complex cellular shapes, the proposed approach models both intracellular and extracellular spaces as distinct poroelastic media, separated by a permeable cell membrane which allows hydrostatic and osmotic pressure-driven fluid exchange. Based on a three-field (displacement, total pressure, and fluid pressure) formulation, the solver leverages the framework of norm-equivalent preconditioning and appropriately fitted norms to ensure robustness across all material parameters of the model. Scalability for large and complex geometries is achieved through efficient Algebraic Multigrid (AMG) approximations of the preconditioners' individual blocks. Furthermore, we accommodate diverse boundary conditions, including full Dirichlet boundary conditions for displacement, which we handle efficiently using the Sherman-Morrison-Woodbury formula.
  Our theoretical analysis is complemented by numerical experiments demonstrating the preconditioners' robustness and performance across various parameters relevant to realistic scenarios. A large scale example of cellular swelling on a dense reconstruction of the mouse visual cortex highlights the method's potential for investigating complex physiological processes such as cellular volume regulation in detailed biological structures.
\end{abstract}

% REQUIRED

% REQUIRED

\section{Introduction}

Mechanical forces are ubiquitous in neuro-biological processes, governing cellular mechanisms such as neuron growth, synapse formation, and brain tissue development \cite{pillai2024mechanics, tyler2012mechanobiology}. These forces, including fluid flow, mechanical stress, and tissue deformation, are crucial to understanding how physical interactions can affect brain function and pathology. For instance, cellular swelling and volume regulation modulate brain states \cite{rasmussen2020interstitial}, while mechanical forces play a pivotal role in structural plasticity and learning \cite{quintana2023mechanical}. Recent advances in high-resolution electron microscopy provide reconstructions of the mammalian brain with unprecedented cellular detail~\cite{turner2022reconstruction}. This paves the way for new numerical investigations into the relationship between complex cellular shape and physiological function. However, performing such simulations on these highly detailed geometries requires efficient and scalable methods for the large-scale linear systems that arise from discretizing the governing equations of cellular mechanics.

\subsection{Governing equations}
\begin{figure}
    \centering
    \includegraphics[width = 0.4 \textwidth]{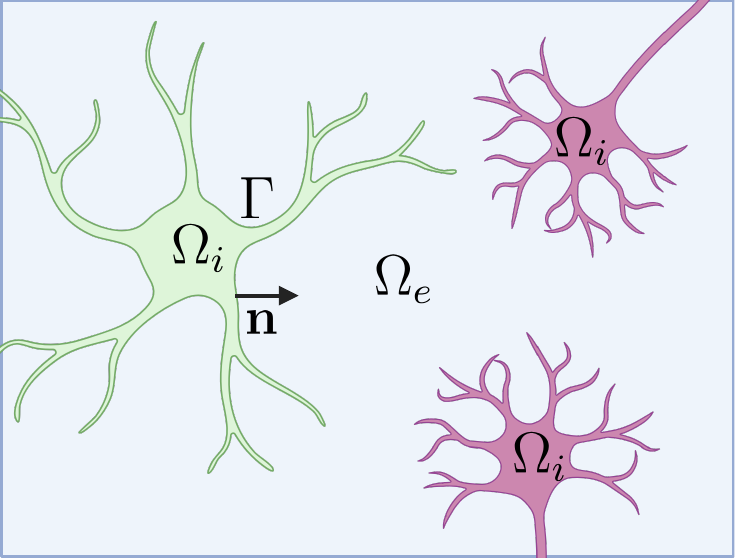}
    \hspace{0.8cm}
    \includegraphics[width = 0.4 \textwidth]{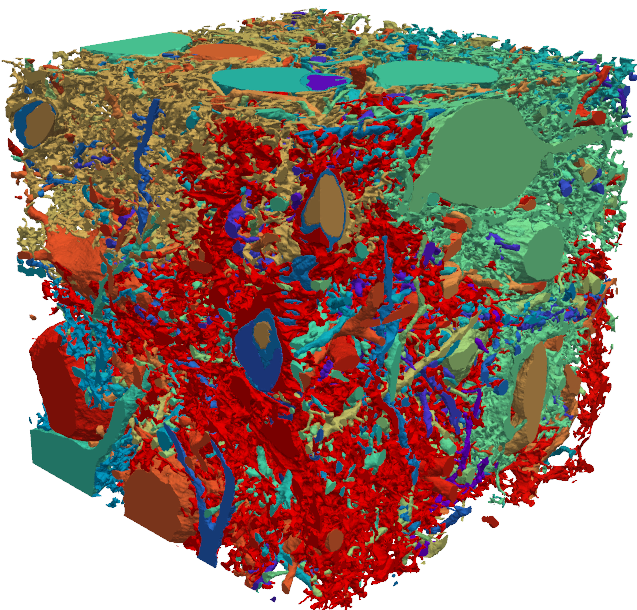}
    \caption{Left) illustration of the geometrical setting: $\Omega_i$ describes the intracellular spaces, separated from the extracellular space $\Omega_e$ by the membrane $\Gamma$. The interface normal $\n$ points from the intra- into the extracellular domain; right) rendering of a realistic tissue reconstruction from the mouse brain, containing 200 cells in a tissue cube with $20\,\mu$m sidelength.}
    \label{fig:geom}
\end{figure}

To model the mechanical interactions between cells and their surrounding extracellular space, we describe the tissue as a composite of distinct poroelastic domains. Specifically, we represent $N$ intracellular regions (physiological cells) by $N$ disjoint domains $\Omega_{i^n} \subset \mathbb{R}^d$ with $d \in \{2,3\}$ for $n=1,\dots, N$ embedded in the extracellular space $\Omega_e \subset \mathbb{R}^d$, and denote the complete intracellular space by $\Omega_i = \bigcup_{i=1,\dots,N} \Omega_{i^n}$ (see \cref{fig:geom}). Further, we define the computational domain by $\Omega = \Omega_i \cup \Omega_e$ with boundary $\partial \Omega$, the cell membrane by $\Gamma = \partial\Omega_i \cap \partial\Omega_e$ and let $\n$ be the unit normal vector pointing from $\Omega_i$ to $\Omega_e$ on $\Gamma$ and outwards on $\partial\Omega$. We characterize both the intracellular space $\Omega_i$ and the extracellular space $\Omega_e$ as deformable porous media consisting of a porous elastic solid filled with interstitial fluid, representing the cytoskeleton or the extracellular matrix, respectively \cite{moeendarbary2013cytoplasm, esteki2021poroelastic}.
Both materials are parameterized by the hydraulic conductivity $\kappa_k > 0$, the storage coefficient $c_{0,k} \geq 0$, the Biot-Willis coefficient $\alpha_k > 0$ and the Lamé constants $\lambda_k > 0$ and $\mu_k > 0$ in their respective domain $\Omega_k$ for $k \in \{i,e\}$. Then, the primary unknown fields are the displacement $\d_k$ and the fluid pressure $p_{F,k}$ defined in each domain $\Omega_k$. 
Simplifying notation, we denote variables (e.g., displacement $\d$, fluid pressure $p_F$) and parameters (e.g., $\kappa, \mu, \lambda, \alpha, c_0$) previously defined with subscript $k \in \{i,e\}$ for domains $\Omega_k$ without the subscript in equations applicable to either domain. These symbols represent global fields, defined by their domain-specific values: $\phi(x) := \phi_k(x) \text{ if } x \in \Omega_k, \text{for } k \in \{i,e\}$. The governing equations are thus written as 
\begin{align}
\label{eq:full_biot}
\begin{aligned}
- \div [ 2 \mu \eps (\d) + \lambda \div \d \I  - \alpha p_F \I] &= \f \quad \text{ in } &\Omega_i \cup \Omega_e, \\
c_0 \partial_t p_F + \alpha \div \partial_t \d - \div (\kappa \nabla p_F)  &= g \quad \text{ in } &\Omega_i \cup \Omega_e.
 \end{aligned}
\end{align}
Here, $\eps(\d)$ is the symmetric gradient of $\d$, $\I$ is the identity matrix, and $\f$ and $g$ denote body forces and fluid sources, respectively.
We denote the jump of a variable $\phi$ across the membrane $\Gamma$ by $\jump{\phi}$, and define it by 
\begin{equation}
    \jump{\phi} = \phi_i|_{\Gamma} -  \phi_e|_{\Gamma}.
\end{equation}
Next, we impose continuity of displacement, conservation of (fluid) mass and momentum as interface conditions on $\Gamma$:
\begin{equation}\label{eq:other_iface_conditions}
\jump{\d} = 0, \quad \jump{\kappa \nabla p_F \cdot \n} = 0, \quad \text{ and } \jump{(2 \mu \eps(\d) + \lambda \div \d - \alpha p_F \mathbf{I})\cdot \n} = 0.
\end{equation}
Allowing for a transmembrane pressure jump and pressure-driven fluid flow across the cell membrane, we impose a Robin-type interface condition for the transmembrane flux $- \kappa \nabla p_F$,
\begin{align}
\label{eq:interface_flow}
- \kappa \nabla p_F\cdot \n = L_p ( \jump{p_F} + p_{\rm osm}) \text{ on }\Gamma,
\end{align}
where $L_p > 0$ is the membrane permeability and $p_{\rm osm}$ is a given osmotic pressure difference. The system is complemented with suitable initial and boundary conditions, specified as needed throughout the remainder of the paper.

The need for efficient iterative solvers for the model is driven by the geometrical complexity of the problem setting. The tissue of the brain features complex structures of intertangled cells with different morphologies, with neuron dendrites and astrocyte processes forming extensive networks (see \Cref{fig:geom}). Representing these structures explicitly leads to computational meshes with large numbers of cells and vertices, requiring scalable, order-optimal methods to solve the large linear systems resulting from discretization of
\eqref{eq:full_biot}-\eqref{eq:interface_flow}.
Furthermore, the model incorporates a wide range of material parameters, which can vary by several orders of magnitude. For instance, in applications related to cellular mechanics in neuroscience, the (dimensionless) effective hydraulic conductivity ($\kappa$) can range from $10^{-7}$ to $10^3$, the (dimensionless) Lamé parameter ($\lambda$) from 10 to $10^5$, and the (dimensionless) membrane permeability ($L_p$) from $10^{-9}$ to $10^{-2}$. An in-depth discussion of the relevant material parameters, for both applications in neuroscience and other fields, will be given in \Cref{sec:parameter_values}. 
Iterative solvers often exhibit sensitivity to such substantial parameter variations, and thus require parameter-robust preconditioning techniques. Developing such techniques for the coupled, multi-domain case is the central goal of this work.

In contrast to our multi-domain poroelasticity model, the development of effective numerical methods for single-domain Biot's
model has been an active area of research, providing a foundation for the present work. This rich body of literature includes
significant efforts focused on stable discretizations and efficient solvers in the form of monolithic preconditioners or
sequential/split approaches (see e.g. \cite{both2017robust, mikelic2013convergence, kim2011stability}). In addition to
conforming discretizations of the classical two-field formulation~\cite{adler2018robust}, previous work includes the development of nonconforming flux
formulations, such as the Crouzeix-Raviart element for displacement coupled with flux condensation techniques \cite{hu2017nonconforming} and discontinuous approximations to ensure global and local mass conservation~\cite{hong2017parameter}. A four-field formulation,
introducing both flux and total pressure as additional unknowns, has been proposed to achieve stability in both the incompressible
and low-permeability limits \cite{boon2021robust}. Other approaches have focused on stabilization techniques, for example, using
a facet bubble function to enrich a first order continuous displacement approximation, and achieving a low number of degrees of
freedom though condensation of bubble degrees of freedom \cite{ohm2019new,rodrigo2018new}. More broadly, various discretization
and preconditioning strategies for Biot's consolidation model continue to be an active topic of discussion and development \cite{rodrigo2023parameter}.

In this paper, we employ a three field formulation with total pressure as additional unknown, as presented in \cite{lee2017parameter}, and extend it to the case of two coupled poroelastic domains. We demonstrate that a straight-forward application of their result to each of the single Biot subdomains is not sufficient to achieve fully parameter-independent stability estimates and
in turn robust preconditioning. Additionally, in contrast to \cite{lee2017parameter}, we refrain from fully diagonalizing
the pressure block, which allows to cover the general case without additional assumptions on the scaling behavior of the
material parameters.

The paper is structured as follows: \Cref{sec:total_pressure} details the total pressure formulation of the cell-by-cell model, including the time discretization and rescaling of the governing equations. \Cref{sec:stability} is dedicated to the stability analysis of the continuous problem, where we establish parameter-independent stability using the framework of fitted norms and operator preconditioning, and discuss the specific case of full Dirichlet boundary conditions for the displacement variable. In \Cref{sec:discr_mg}, we describe the finite element discretization of the model and the construction of efficient algebraic multigrid (AMG) approximations for the preconditioners' individual blocks. \Cref{sec:numerical_results} presents numerical results, including 2D and 3D examples demonstrating parameter robustness, the performance of multigrid approximations, and a large-scale simulation of cellular swelling in a realistic brain geometry. 

\section{Total pressure formulation of the poroelastic cell-by-cell model}\label{sec:total_pressure}

We begin our analysis by simplifying the time-dependent system \eqref{eq:full_biot} into a static problem at a single time step. This is achieved through time discretization and appropriate rescaling of the parameters.
Employing the backward Euler method with time step $\tau > 0$, the second equation of \eqref{eq:full_biot} at time $t_k = t_{k-1} + \tau$ becomes:
\begin{equation}
        \frac{c_0}{\tau} (p_F^k - p_F^{k-1}) + \frac{\alpha}{\tau}\left(\div \d^k - \div \d^{k-1} \right) - \div(\kappa \nabla p_F^k) = g^k. 
\end{equation}
For clarity in the subsequent stability analysis, we rescale the entire system \eqref{eq:full_biot} by $2\mu$ and group all terms from previous time steps into a single source term $\tilde{g}$. This leads to the static, rescaled momentum and mass balance equations
\begin{equation}
\begin{aligned}
\label{eq:strong}
- \div \eps (\d) - \tilde{\lambda} \nabla \div \d  + \tilde{\alpha}\nabla p_F &= \tilde{\f} \quad \text{ in } &\Omega,\\
 - \tilde{c_0} p_F - \tilde{\alpha} \div \d + \div (\tilde{\kappa} \nabla p_F) &=  \tilde{g} \quad \text{ in } &\Omega,\\
\end{aligned}
\end{equation}
while the Robin interface condition becomes
\begin{equation}\label{eq:strong_iface}
  - \tilde{\kappa} \nabla p_F \cdot \n = \tilde{L}_p ( \jump{p_F} + p_{\rm osm})\quad\text{ on }\Gamma.
\end{equation}
Here the rescaled parameters are 
\begin{equation}
\tilde{\lambda} = \frac{\lambda}{2 \mu}, \quad
\tilde{\alpha} = \frac{\alpha}{2 \mu}, \quad
\tilde{\kappa} = \frac{\kappa \tau}{2 \mu}, \quad
\tilde{c_0} = \frac{c_0}{2 \mu}, \quad
\tilde{L_p} = \frac{L_p \tau}{2 \mu}, \quad
\tilde{\f} = \frac{\f}{2 \mu}, \quad
\end{equation}
and $\tilde{g}$ incorporates all terms from $t_{k-1}$. To simplify notation, we drop the tilde symbol for the remainder of the paper.

\subsection{Total pressure formulation}

Motivated by its stability in the incompressible limit of the solid matrix ~\cite{lee2017parameter, lee_parameter-robust_2016}, we employ a total pressure formulation of Biot's consolidation model with the additional scalar unknown $p_T = - \lambda \div \d + \alpha p_F$. We can thus rewrite system \eqref{eq:strong}-\eqref{eq:strong_iface} as follows
\begin{equation}
\begin{aligned}\label{eq:TP_strong}
- \div \eps (\d)  + \nabla p_T &= \f \quad \text{ in } &\Omega,\\
- \div \d - \lambda^{-1} (p_T - \alpha p_F) &= 0 \quad \text{ in } &\Omega,\\
\lambda^{-1} \alpha p_T - (\alpha^2 \lambda^{-1} + c_0) p_F + \div (\kappa \nabla p_F) &= g \quad \text{ in } &\Omega,\\
- \kappa \nabla p_F^e \cdot \n  = - \kappa \nabla p_F^i \cdot \n &=L_p ( \jump{p_F} + p_{\rm osm}) \quad \text{ on } &\Gamma.
\end{aligned}
\end{equation}
We equip the system \eqref{eq:TP_strong} together with \eqref{eq:other_iface_conditions} with suitable boundary conditions. Assuming two partitions of the boundary $\partial \Omega$,
\begin{equation*}
\partial \Omega = \Gamma_p \cup \Gamma_f, \quad \partial \Omega = \Gamma_d \cup \Gamma_t,
\end{equation*}
we impose the boundary conditions
\begin{align*}
p_F &= 0 \text{ on } \Gamma_p, &  -\kappa \nabla p_F \cdot  \n &= 0 \text{ on } \Gamma_f, \\
\d &= \mathbf{0} \text{ on } \Gamma_d, & \left[ \eps(\d) - p_T  \I \right] \cdot \n &= \mathbf{0} \text{ on } \Gamma_t.
\end{align*}
We require that $|\Gamma_d| > 0$ and $|\Gamma_t| \geq 0$, i.e. that a Dirichlet condition for the displacement is applied on a nonzero part of the boundary, while the traction boundary condition is optional. This requirement eliminates rigid body motions which would otherwise render the elasticity operator singular. While the pure traction case ($| \Gamma_d |=0$) can be treated by removing the rigid body modes via Lagrange multipliers or using orthogonalization techniques to resolve the kernel of the discrete operator (see, e.g., \cite{kuchta2019singular}), we exclude it here to focus the analysis on the stability of the coupled multi-domain problem. Additionally, we allow both $|\Gamma_p|=0$ or $|\Gamma_f| = 0$, i.e. pure Dirichlet or Neumann boundary conditions for the fluid pressure.

Based on the standard function spaces $H^1(\Omega)$, $\mathbf{H}^1(\Omega) = H^1(\Omega)^d$, and $L^2(\Omega)$, we define the spaces
\begin{align*}
    \mathbf{V} = \mathbf{H}^1_{0, \Gamma_d}(\Omega), \quad Q_T = L^2(\Omega), \quad Q_F = \{q_F=(q^i_F, q^e_F): q^k_F\in H^1_{0, \Gamma_p}(\Omega^k)) \}
\end{align*}
reflecting the continuity of displacement on $\Gamma$, the essential boundary conditions and the discontinuous fluid pressure across the interface.

Further, let $(\cdot,\cdot)$ and $(\cdot,\cdot)_{\Gamma}$ denote the $L^2$ inner product of scalar, vector or matrix valued functions on $\Omega$ and $\Gamma$, respectively. 
We can then state the variational formulation of the problem: Find $\d \in \mathbf{V}$, $p_T \in Q_T$ and $p_F := (p_F^e, p_F^i) \in Q_F$ such that
\begin{align}
\label{eq:total_pressure}
\begin{aligned}
(\eps (\d), \eps (\v)) - (\div \v, p_T) &= (\f, \v) \quad &\forall \v \in \mathbf{V}, \\
- (\div \mathbf{d}, q_T) - (\lambda^{-1} p_T, q_T) + (\alpha \lambda^{-1} p_F, q_T)  &= 0 \quad &\forall q_T \in Q_T, \\
(\alpha \lambda^{-1} p_T, q_F) - (\kappa \nabla p_F, \nabla q_F) &- ((\alpha \lambda^{-1} + c_0) p_F, q_F)  \\ - (L_p \jump{p_F}, \jump{q_F})_{\Gamma} &= (g, q_F) + (L_p p_{\rm osm}, \jump{q_F})_{\Gamma} &\forall q_F \in Q_F,\\
\end{aligned}
\end{align}
with $q_F = (q_F^e, q_F^i)$.
Note that the interface term $(L_p \jump{p_F}, \jump{q_F})_{\Gamma}$ stems from integration by parts and application of the interface condition \eqref{eq:interface_flow}:
\begin{equation*}
\begin{aligned}
\int_{\Omega_e \cup \Omega_i} q_F \div (\kappa \nabla p_F) dx
&= - \int_{\Omega_e \cup \Omega_i} \nabla q_F \cdot (\kappa \nabla p_F) dx
+ \int_{\Gamma} q_F^i (\kappa \nabla p_F^i) \cdot \n ds
- \int_{\Gamma} q_F^e (\kappa \nabla p_F^e) \cdot \n ds \\
&= - \int_{\Omega_e \cup \Omega_i} \nabla q_F \cdot (\kappa \nabla p_F)
+ \int_{\Gamma} (q_F^i - q_F^e) \kappa \nabla p_F^i \cdot \n ds  \\
&= - (\kappa \nabla p_F, \nabla q_F) - (L_p\jump{p_F}, \jump{q_F})_{\Gamma} - (L_p p_{\rm osm}, \jump{q_F})_{\Gamma}.
\end{aligned}
\end{equation*}

Defining the interface trace operator $(Tp_F,q_F) := (\jump{p_F}, \jump{q_F})_{\Gamma}$, we can restate the problem \eqref{eq:TP_strong}
in matrix/operator notation as
\begin{align}\label{eq:matrix_form}
\mathcal{A} 
\begin{bmatrix}
\d \\
p_T \\
p_F
\end{bmatrix} :=
\begin{bNiceArray}{cIcc}[margin]
- \div \eps  & -\nabla & 0  \\
\dashedline
- \div & - \lambda^{-1}I & \alpha \lambda^{-1}I \\
0 &   \alpha\lambda^{-1} I & - ( \alpha^2 \lambda^{-1} + c_0)I + \div(\kappa \nabla) - L_p T
\end{bNiceArray}
\begin{bmatrix}
\d \\
p_T \\
p_F
\end{bmatrix} = 
\begin{bmatrix}
\f \\
0 \\
g
\end{bmatrix}.
\end{align}

\begin{example}[Na\"ive preconditioner]\label{sec:motivating_example}
To highlight the necessity of a tailored preconditioner for the coupled problem, we first consider a "na\"ive" approach where we apply the block-diagonal preconditioner suggested by Lee et al by Lee et al~\cite{lee2017parameter} for a single Biot domain to the coupled system \eqref{eq:matrix_form}. The preconditioner given by
\begin{equation}\label{eq:prec_naiv}
\mathcal{B}_{single} := 
\begin{bmatrix}
(- \div \eps)^{-1} & 0 & 0\\
0 & I^{-1} & 0 \\
0 & 0&  (\alpha^2 \lambda^{-1}I - \div \kappa \nabla)^{-1}
\end{bmatrix}
\end{equation}
naturally lacks the interface coupling terms involving $L_p$ present in the coupled operator. Exploring its stability in the coupling parameter $L_p$ on a unit square geometry with the interface $\Gamma$ at $x=0.5$, we find that the method is strongly sensitive to variations  in $L_p$ (see ~\Cref{fig:motivating_example}). This sensitivity underscores the need to account for the interface term when designing a parameter-robust preconditioner for the coupled system, which is the focus of the next section.

\begin{figure}
    \centering
    \includegraphics[width=0.8\linewidth]{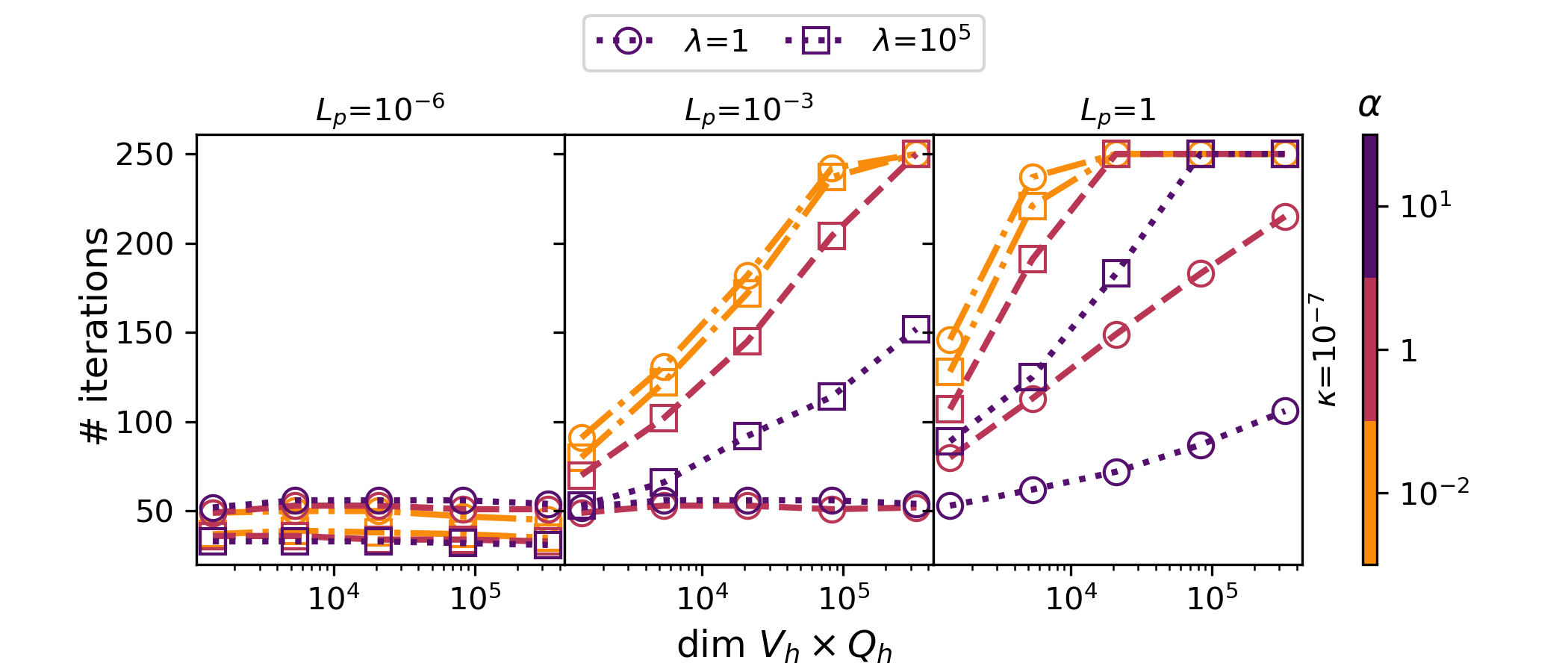}
    \caption{Number of MinRes iterations of system \eqref{eq:total_pressure} with preconditioner \eqref{eq:prec_naiv} on a unit square geometry. The  number of iterations is capped at 250.}
    \label{fig:motivating_example}
\end{figure}
\end{example}

\section{Parameter independent stability of the continuous poroelastic cell-by-cell}\label{sec:stability}

This section focuses on establishing the parameter-independent stability of the continuous poroelastic cell-by-cell model.
We begin by revisiting the theoretical framework of operator preconditioning and the concept of ``fitted'' norms, which are
crucial for analyzing the stability of perturbed saddle-point problems. Subsequently, we apply the framework to our
specific cell-by-cell model and establish its uniform stability in suitably weighted norms. 

\subsection{Operator preconditioning}

We briefly recall the framework of parameter-robust operator preconditioner as presented in \cite{mardal_preconditioning_2011}. Let $X$ be a separable, real Hilbert space with inner product $(\cdot, \cdot)_X$, and denote its dual by $X^{*}$. Further, let $L(X,X^{*})$ be the space of bounded linear operators mapping from $X$ to $X^{*}$. For an invertible and symmetric operator $\mathcal{A} \in L(X,X^*)$ and a given $\mathcal{F} \in X^*$, consider the problem: Find $x \in X$ such that
\begin{equation}
    \mathcal{A} x = \mathcal{F}.
\end{equation}
Introducing the symmetric isomorphism $\mathcal{B} \in L(X^*, X)$, the preconditioned problem becomes:
\begin{equation}
    \mathcal{B} \mathcal{A} x = \mathcal{B} \mathcal{F}.
\end{equation}
The convergence rate of a Krylov-subspace method applied to this problem is bounded by the condition number $K(\mathcal{BA})$, which is defined as 
\begin{equation}
    K(\mathcal{BA}) = ||\mathcal{BA}||_{L(X,X)} || (\mathcal{BA})^{-1}||_{L(X,X)}.
\end{equation}
For a parameter dependent operator $\mathcal{A}_\epsilon \in L(X_\epsilon, X^*_\epsilon)$, we aim to find a preconditioner operator $\mathcal{B}_\epsilon \in L(X^*_\epsilon, X_\epsilon)$ such that the condition number is uniformly bounded in a set of parameters $\epsilon$. We can systematically construct such a parameter robust preconditioner by choosing the inner products on the weighted Hilbert spaces $X_\epsilon$ and $X^*_\epsilon$ such that $||\mathcal{A}_\epsilon||_{L(X_\epsilon, X^*_\epsilon)}$ and $||\mathcal{B}_\epsilon||_{L(X^*_\epsilon, X_\epsilon)}$
are bounded independently of $\epsilon$. Then, a canonical preconditioner $B_\epsilon$ corresponds to the Riesz map induced by
the inner product, and the condition number $K(\mathcal{B}_\epsilon \mathcal{A}_\epsilon)$ is uniformly bounded in $\epsilon$.
We next review how to establish $\mathcal{A}$ as isomorphism in our specific context when to operator arises due to
a perturbed saddle point problem.

\subsection{Stability analysis of perturbed saddle-point problems}

For the stability analysis of perturbed saddle point problems, we apply the framework recently proposed in
\cite{hong2023new} (see also \cite{braess1996stability}). For convenience, we restate the main results here.
We assume $V$ and $Q$ are two Hilbert spaces and set $X:= V \times Q$. The framework establishes stability conditions for saddle-point problems induced by bilinear forms $\mathcal{A}(\cdot, \cdot)$ on $X \times X$ of the form
\begin{equation}
\label{eq:abstract_sp}
\mathcal{A}((u, p), (v,q)) = a(u,v) + b(v,p) + b(u,q) - c (p,q),
\end{equation}
with symmetric positive semi-definite forms $a$ on $X$, and $c$ on $Q$ and a bilinear form $b(\cdot, \cdot)$ on $V \times Q$.
For $\mathcal{F} \in X'$, the abstract saddle-point problem then reads: Find $(u, p) \in V \times Q$ such that
\begin{equation}
    \mathcal{A}((u, p), (v,q)) = a(u,v) + b(v,p) + b(u,q) -c (p,q) =  \mathcal{F}((v,q)) \quad \forall v \in V, \forall q \in Q.
\end{equation}
A key aspect of the framework is the construction of appropriate ``fitted'' norms, which are defined in the following: 
\begin{definition}[Fitted norms \cite{hong2023new}]
\label{def:fitted_norms}
For two Hilbert spaces $V$ and $Q$, the norms $|| \cdot ||_V$ on $V$ and $|| \cdot ||_Q$ on Q are called fitted if they satisfy the splittings
\begin{align}
    ||q||_Q^2 &= |q|^2_Q + c(q,q) =: \langle \hat Q q,q \rangle_{Q' \times Q},   \\
    ||v||_V^2 &= |v|^2_V + |v|_b^2,
\end{align}
where $|\cdot|_V$ and $|\cdot|_Q$ are seminorms on $V$ and $Q$, $\hat{Q}$ is a linear operator $\hat{Q}:Q' \rightarrow Q$ defined by $|| \cdot ||_Q$ and $|\cdot|_b$ is a seminorm on V with 
\begin{equation}
    |v|_b^2 = \langle B v, \hat Q^{-1} B v\rangle_{Q' \times Q }.
\end{equation}
\end{definition}
Then, the main stability result (Th. 2.14, \cite{hong2023new}) follows.

\begin{theorem}[\cite{hong2023new}, Th. 2.14]
\label{th:abstract_saddle_point}
Let $V$ and $Q$ be Hilbert spaces with the fitted norms $|| \cdot ||_V$ and $|| \cdot ||_Q$ and let $\mathcal{A}$  be the linear operator induced by \eqref{eq:abstract_sp}. Assume further that $a(\cdot, \cdot)$ is continuous, $a(\cdot, \cdot)$ and $c(\cdot, \cdot)$ are symmetric positive semi-definite. The continuity of $b(\cdot, \cdot)$ and $c(\cdot, \cdot)$ in the fitted norms directly follows from their definitions. If $a(\cdot, \cdot)$ satisfies the coercivity estimate 

\begin{equation}
\label{eq:coercitvity}
    a(v, v) \ge C_a |v|^2_V \quad \forall v \in V
\end{equation}
and there exists a constant $\beta > 0$ such that
\begin{equation}
\label{eq:small_infsup}
    \sup_{\substack{v \in V \\ v \neq 0}} \frac{b(v,q)}{||v||_V} \ge \beta |q|_Q \quad \forall q \in Q
\end{equation}
then, the bilinear form $\mathcal{A}(\cdot, \cdot)$ is continuous and inf-sup stable under the combined norm
$||\cdot||_X$, $||\cdot||^2_X = ||\cdot||^2_V + ||\cdot||^2_Q$.
    
\end{theorem}

We use this result and the idea of fitted norms to prove the parameter-independent stability of the model and guide the design of parameter-robust preconditioners.

\subsection{Stability of the cell-by-cell model with mixed displacement boundary conditions}

We begin with the case where both $|\Gamma_d| > 0$ and $|\Gamma_t| > 0$, i.e. when both a Dirichlet condition for the displacement and a traction boundary condition are applied on a nonzero part of the boundary.
Casting the total pressure formulation of the poroelastic cell-by-cell model into the abstract
saddle point framework, the bilinear forms defining the operator $\mathcal{A}$ in \eqref{eq:abstract_sp} read:
\begin{align}
    a(\d, \v) &= (\eps (\d), \eps (\v)) \quad &\forall \d, \v \in \mathbf{V}, \\
    b(\v, q) &= -(\div\v, q_T) \quad &\forall \v \in \mathbf{V}, \forall q \in Q, \\
    c(p, q) &= (\lambda^{-1} p_T, q_T) -(\alpha \lambda^{-1} p_T, q_F) - (\alpha \lambda^{-1} p_F, q_T) + (\kappa \nabla p_F, \nabla q_F) \\
    &+ ((\alpha^2 \lambda^{-1} + c_0) p_F, q_F) + (L_p \jump{p_F}, \jump{q_F})_{\Gamma} \quad &\forall p,q \in Q,  \nonumber
\end{align}
where $Q = Q_T \times Q_F$ with $p=(p_T, p_F)$ and $q=(q_T, q_F)$. 

Next, we turn to suitably fitted norms and start with the semi-norms:
\begin{align}
 |q|^2_Q &:= (q_T, q_T) \quad &\forall q \in Q, \\
 |\v|^2_V &:= (\eps (\v), \eps (\v)) \quad &\forall \v \in \mathbf{V}.
\end{align}
This results in the following full norms on $\mathbf{V}$ and $Q$:
\begin{equation}\label{eq:tp_norms}
\begin{aligned}
  ||q||^2_Q &= |q|^2_Q +  c(q,q) = (q_T, q_T) + (\lambda^{-1} q_T, q_T) -(\alpha\lambda^{-1} q_T, q_F) - (\alpha\lambda^{-1} q_F, q_T)\\
  & + (\kappa \nabla q_F, \nabla q_F)
    + ((\alpha^2\lambda^{-1} + c_0) q_F, q_F) + (L_p \jump{q_F}, \jump{q_F})_{\Gamma} \quad \forall q \in Q, \\
||\v||^2_V &= |\v|^2_V + \langle B \v, \hat{Q}^{-1} B \v \rangle _{Q' \times Q} = (\eps (\v), \eps (\v)) + \langle B \v, \hat{Q}^{-1} B \v \rangle _{Q' \times Q} \quad \forall \v \in \mathbf{V},
\end{aligned}
\end{equation}
where $B:V \rightarrow Q'$ is the operator $B:= \begin{bmatrix}
    - \div \\
    0
\end{bmatrix}$.

\begin{theorem}[Parameter robust stability of the continuous total pressure formulation]\label{th:total_pressure_stability}

The saddle point problem \eqref{eq:total_pressure} is continuous and inf-sup stable independent of the material parameters $\alpha$, $\kappa$, $\lambda$, $c_0$ and $L_p$ under the combined norm $|| \cdot ||_X = (|| \cdot ||^2_V + || \cdot ||^2_Q)^{1/2}$ with $|| \cdot ||_V$ and $|| \cdot ||_Q$ as defined in \eqref{eq:tp_norms}.
    
\end{theorem}

\begin{proof}
Applying \Cref{th:abstract_saddle_point}, we first verify continuity of the symmetric bilinear form $a(\cdot, \cdot)$ by observing that
\begin{equation*}
|a(\v, \u)| = |(\eps(\v), \eps(\u))| \leq |\v|_V |\u|_V \leq ||\v||_V ||\u||_V \quad \forall \u, \v \in \mathbf{V}
\end{equation*}
and note that the coercivity estimate \eqref{eq:coercitvity} holds
\begin{equation*}
 a(\v, \v) = |\v|_V^2 \quad \forall \v \in \mathbf{V},
\end{equation*}
which also implies positive-definiteness. Further, $c(\cdot, \cdot)$ is clearly symmetric, and we further observe that
\begin{equation*}
 (\lambda^{-1} p_T, p_T) -2(\alpha \lambda^{-1} p_T, p_F) + (\alpha^2\lambda^{-1} p_F, p_F) = 
(\lambda^{-1}(p_T - \alpha p_F), p_T - \alpha p_F) \geq 0
\end{equation*}
and thus 
\begin{align}\label{eq:c_pos_semipef}
  c(p, p) &= (\lambda^{-1} p_T, p_T) -2(\alpha \lambda^{-1} p_T, p_F) + ((\alpha^2\lambda^{-1} + c_0) p_F, p_F) + (\kappa \nabla p_F, \nabla p_F) + \\
  & (L_p \jump{p_F}, \jump{p_F})_{\Gamma} \geq (c_0 p_F, p_F) + (\kappa \nabla p_F, \nabla p_F)
+ (L_p \jump{p_F}, \jump{p_F})_{\Gamma}
\geq 0 \quad \forall p \in Q, \notag
\end{align}
which establishes positive semi-definiteness of $c(\cdot, \cdot)$.

It remains to show the inf-sup condition \eqref{eq:small_infsup}. To this end we construct $\v_0 \in \mathbf{V}$ such
that $- \div \v_0 = q_T$ and for which it also holds that $||\v_0||_{H^1} \leq \gamma || q_T||_{L^2}$ for some $\gamma > 0$.
Note that the existence of such a $\v_0$ follows from the classical inf-sup condition for the Stokes problem with outflow boundary condition, see~\cite{braack2025inf}. We find:
\begin{equation}\label{eq:stokes_est}
b(\v_0, q) = (q_T, q_T) = |q|^2_Q
\text{ and } 
||\v_0||_{H^1} \leq \gamma ||q_T ||_{L^2} = \gamma |q|_Q.
\end{equation}
We also observe that:
\begin{equation*}
\begin{aligned}
 \langle B \v, \hat{Q}^{-1} B \v \rangle _{Q' \times Q} 
 &= ||B\v||^2_{Q'} 
 = \left( \sup_{q \in Q} \frac{b(\v, q)}{||q||_Q} \right)^2
 = \left( \sup_{q \in Q} \frac{(\div \v, q_T)}{||q||_Q} \right)^2
 \leq \left( \sup_{q \in Q} \frac{||\div \v|| || q_T ||}{||q||_Q} \right)^2 \\
 &\leq (\div \v, \div \v) \quad \forall \v \in V.
\end{aligned}
\end{equation*}
In turn,
\begin{equation}\label{eq:V_bound}
\begin{aligned}
||\v_0 ||^2_V &= (\eps (\v_0), \eps (\v_0)) + \langle B \v_0, \hat{Q}^{-1} B \v_0 \rangle _{Q' \times Q}
\leq (\eps (\v_0), \eps (\v_0)) + (\div \v_0, \div \v_0) \\
&\leq 2 ||\v_0||^2_1 
\leq 2 \gamma^2 |q|_Q^2.
\end{aligned}
\end{equation}
Finally, we verify condition \eqref{eq:small_infsup}
\begin{equation}
\sup_{\v \in V} \frac{b(v,q)}{||\v||_V}
\geq \frac{b(\v_0, q)}{||\v_0||_V}
\geq \frac{|q|_Q^2}{\sqrt{2} \gamma |q|_Q}
=: \beta |q|_Q
\end{equation}
with $\beta >0 $ independent of the material parameters, which concludes the proof. 
\end{proof}

\Cref{th:total_pressure_stability} implies that the operator $\mathcal{A}$ is invertible and the operator norm is independent of the five parameters $c_0$, $\lambda$, $\kappa$, $\alpha$ and $L_p$. As a direct consequence of this stability, the norm-equivalent preconditioner is parameter-robust.

We further note that $||\cdot||_V$ is equivalent to $|\cdot|_V$ by \eqref{eq:V_bound}. Hence, a norm-equivalent preconditioner is given by
\begin{equation}\label{eq:TP_preconditioner}
\mathcal{B} := 
\begin{bmatrix}
(- \div \eps)^{-1} & \\
& \begin{pmatrix}
    \lambda^{-1}I + I & -\alpha \lambda^{-1}I \\
    -\alpha \lambda^{-1}I & (\alpha^2 \lambda^{-1} + c_0)I - \div \kappa \nabla + L_p T
\end{pmatrix}^{-1}
\end{bmatrix}.
\end{equation}

\begin{remark}\label{rem:diag}

It is possible to diagonalize the preconditioner \eqref{eq:TP_preconditioner} by noting that by \eqref{eq:c_pos_semipef}
\begin{equation*}
||p||_Q^2 = |p|_Q^2 + c(p,p) \geq (p_T, p_T) + (c_0 p_F, p_F) + (\kappa \nabla p_F, \nabla p_F) + (L_p\jump{p_F}, \jump{p_F}) \quad \forall p \in Q,
\end{equation*}
which yields the preconditioner
\begin{equation}\label{eq:TP_preconditioner_diag}
\mathcal{B_{\rm diag}} := 
\begin{bmatrix}
(- \div \eps)^{-1} & \\
& I^{-1} & \\
&  & ( c_0 I - \div \kappa \nabla + L_p T)^{-1}
\end{bmatrix}.
\end{equation}
However, this is generally not possible for purely essential boundary conditions on the displacement variable (see \Cref{rem:diag-P0}). For the sake of a unified presentation of the two scenarios, we therefore focus on the
preconditioner \eqref{eq:TP_preconditioner} for the remainder of the paper.
\end{remark}

\subsection{Stability of the cell-by-cell model with full Dirichlet displacement boundary conditions}\label{sec:full_disp_bc}

In case of Dirichlet boundary conditions for the displacement on the entire domain boundary ($\partial\Omega = \Gamma_d$), we need to slightly modify the fitted norms and thus the resulting preconditioner. In particular, we chose the $Q$-seminorm by
\begin{align}
 |q|^2_Q &:= (q_{T,0}, q_{T,0}) \quad &\forall q \in Q,
\end{align}
where $q_{T,0}=P_0 q_T$ is the $L^2$-projection of $q_T \in L^2(\Omega)$ to $L^2_0(\Omega) = \{ q \in L^2(\Omega): \int_{\Omega}q = 0 \}$.
The choice yields the following full norms on $Q$ and $V$:
\begin{equation}\label{eq:tp_norms_P0}
\begin{aligned}
||q||^2_Q &= (q_{T,0}, q_{T,0}) + (\lambda^{-1} q_T, q_T) -(\alpha \lambda^{-1} q_T, q_F) - (\alpha \lambda^{-1} q_F, q_T) + (\kappa \nabla q_F, \nabla q_F) \\
&+ ((\alpha^2 \lambda^{-1} + c_0) q_F, q_F) + (L_p \jump{q_F}, \jump{q_F})_{\Gamma}  \quad \forall q \in Q,\\
||\v||^2_V &= (\eps (\v), \eps (\v)) + \langle B \v, \hat{Q}^{-1} B \v \rangle _{Q' \times Q} \quad \forall \v \in \mathbf{V}.
\end{aligned}
\end{equation}

\begin{theorem}[Parameter robust stability - homogeneous Dirichlet displacement boundary conditions]\label{thm:D}
The saddle point problem \eqref{eq:total_pressure} equipped with homogeneous Dirichlet boundary conditions on the displacement ($\Gamma_d=\partial\Omega$) is continuous and inf-sup stable independent of the material parameters $\alpha$, $\kappa$, $\lambda$, $c_0$ and $L_p$ under the combined norm $|| \cdot ||_X = (|| \cdot ||^2_V + || \cdot ||^2_Q)^{1/2}$ with $|| \cdot ||_V$ and $|| \cdot ||_Q$ as defined in \eqref{eq:tp_norms_P0}.
\end{theorem}

\begin{proof}
The proof is analogous to the proof of  \Cref{th:total_pressure_stability}, except for \eqref{eq:stokes_est}.
We establish a similar result by choosing $\v_0 \in V$ such that $- \div \v_0 = q_{T,0}$, $||\v_0||_{H^1} \leq \gamma || q_T||_{L^2}$ (see \cite{girault1986finite}, Th. 5.1), which yields
\begin{equation}
b(\v_0, q) = (q_{T,0}, q_{T,0}) = |q|^2_Q
\text{ and } 
||\v_0||_1 \leq \gamma ||q_{T,0} || = \gamma |q|_Q 
\end{equation}
and thus the inf-sup condition is satisfied.
The remaining conditions can be verified analogously to the proof of \Cref{th:total_pressure_stability}.
\end{proof}

As a direct consequence of \Cref{thm:D}, the norms defined in equation \eqref{eq:tp_norms_P0} yield the parameter-robust preconditioner
\begin{equation}\label{eq:TP_preconditioner_L20proj}
\mathcal{B}_d := 
\begin{bmatrix}
(- \div \eps)^{-1} & \\
& \begin{pmatrix}
    \lambda^{-1}I + P_0 & -\alpha \lambda^{-1}I \\
    -\alpha\lambda^{-1}I & (\alpha^2 \lambda^{-1} + c_0)I - \div \kappa \nabla + L_p T
\end{pmatrix}^{-1}
\end{bmatrix}
\end{equation}
in case of homogeneous Dirichlet boundary conditions for the displacement variable.

\begin{remark}\label{rem:diag-P0}

With the seminorm $|p|_Q^2=c(p_{T,0}, p_{T,0})$ the estimate in \Cref{rem:diag} becomes
\begin{equation*}
||p||_Q^2 = |p|_Q^2 + c(p,p) \geq (p_{T,0}, p_{T,0}) + (\alpha p_F, p_F) + (\kappa \nabla p_F, \nabla p_F) + (L_p\jump{p_F}, \jump{p_F}) \quad \forall p \in Q,
\end{equation*}
which no longer yields a full norm on $Q$. In this situation, a fully diagonal preconditioner can still be established, but requires the additional assumption that $c_0$ scales as $c_0 \approx \alpha^2 \lambda^{-1}$ (as presented in \cite{lee_mixed_2019}). Then, one can sharpen estimate \eqref{eq:c_pos_semipef} and show that
\begin{align}
c(p, p) \geq \frac{1}{4}(\lambda^{-1} p_T, p_T) + \frac{2}{3}(\alpha^2 \lambda^{-1} p_F, p_F) + (\kappa \nabla p_F, \nabla p_F)
+ (L_p \jump{p_F}, \jump{p_F})_{\Gamma} \quad
 \forall p \in Q, \notag
\end{align}
which implies the diagonal preconditioner:

\begin{equation}\label{eq:TP_preconditioner_diag_P0}
\mathcal{B_{\rm diag-P_0}} := 
\begin{bmatrix}
(- \div \eps)^{-1} & \\
& (\lambda^{-1}I + P_0)^{-1} & \\
&  & ( \alpha^2 \lambda^{-1} I - \div \kappa \nabla + L_p T)^{-1}
\end{bmatrix}.
\end{equation}
\end{remark}

\begin{remark}[Stability in the limit of incompressible constituents]
\label{rem:incompressible_stability}
In the limit of incompressible constituents ($c_0 \to 0$) with pure Neumann fluid boundary conditions ($|\Gamma_p|=0$), the stability of the pressure depends on the displacement boundary conditions.

If traction boundary conditions are present ($|\Gamma_t| > 0$),  the $Q$-norm for the fluid subproblem (see \eqref{eq:tp_norms}) satisfies $||p||_Q^2 \geq ||p_T||^2 + c_0 ||p_F||^2 + \lambda^{-1} ||p_T - \alpha p_F ||^2$. Consequently, the $Q$-norm remains a full norm even when $c_0=0$: the term $||p_T||^2$ controls the total pressure, and the coupling term $\lambda^{-1}\|p_T - \alpha p_F\|^2$ subsequently controls the fluid pressure. Thus, the problem remains well-posed, and the conditions for theorem \ref{th:abstract_saddle_point}
remain satisfied.

However, if the displacement is fully constrained by Dirichlet boundary conditions ($|\Gamma_t|=0$) and pure Neumann conditions are applied to the fluid ($|\Gamma_p|=0$), the fluid pressure is determined only up to an additive constant and $\mathcal{A}$ becomes singular. Specifically, the system admits a non-trivial kernel spanned by the hydrostatic mode $(\d, p_T, p_F) = (0, \alpha, 1)$.
This is reflected in the $Q$-norm: The relevant semi-norm $|p|^2_Q = (p_{T,0}, p_{T,0})$ controls only the deviatoric part of the total pressure; it fails to control the mean value. Thus, the coupling term cannot constrain the mean of the fluid pressure, leaving it determined only up to a constant.
To ensure the well-posedness of the problem in this limit, it is necessary to restrict the fluid pressure to $L^2_0(\Omega)$. In the multidomain setting, this restriction is global; the membrane coupling terms stabilize the intracellular pressures relative to the extracellular space, meaning that fixing the global pressure level (e.g., in the ECS) is sufficient to determine the pressure in all cells without requiring cell-wise boundary conditions. While a detailed analysis and numerical experiments for this specific singular limit are outside the scope of the present work, the treatment of the hydrostatic mode via global constraints is a standard technique in poroelasticity.
\end{remark}

\section{Discretization and multigrid solver}\label{sec:discr_mg}

In this section, we describe the discretization of the problem and present the solvers used for the resulting linear systems. In particular, we choose appropriate multigrid methods for approximate solutions of each of the blocks of the norm-equivalent preconditioner with a focus on efficiency, scalability and availability in the widely adopted PETSc framework \cite{balay2024petsc}.

\subsection{Discretization of the cell-by-cell model}
We discretize system \eqref{eq:total_pressure} with the conforming finite element spaces $\mathbf{V}_h \subset \mathbf{V}$, $Q_{T,h} \subset Q_T$ and $Q_{F,h} \subset Q_F$. The discrete counterpart of \eqref{eq:total_pressure} then reads: Find $(\dh, \pth, \pfh) \in \mathbf{V}_h \times Q_{T,h} \times Q_{F,h}$ such that
\begin{align}
\label{eq:total_pressure:discr}
\begin{aligned}
(\eps (\dh), \eps (\v)) - (\div \v, \pth) &= (\f, \v) \quad &\forall \v \in \mathbf{V}_h, \\
- (\div \dh, q_T) - (\lambda^{-1} p_{T,h}, q_T) + (\alpha \lambda^{-1} p_{F,h}, q_T)  &= 0 \quad &\forall q_T \in Q_{T,h}, \\
(\alpha \lambda^{-1} p_{T,h}, q_F) - (\kappa \nabla p_{F,h}, \nabla q_F) &- ((\alpha \lambda^{-1} + c_0) p_{F,h}, q_F) \\ - (L_p \jump{p_{F,h}}, \jump{q_F})_{\Gamma} = (g, q_F)  &+ (L_p p_{\rm osm}, \jump{q_F})_{\Gamma} &\forall q_F \in Q_{F,h}.\\
\end{aligned}
\end{align}

The stability of the discretization requires the pair $\mathbf{V}_h \times Q_{T,h}$ to be Stokes stable, i.e. 
\begin{equation}
    \inf_{p_{T} \in Q_{T,h}} \sup_{\d \in {\mathbf{V}_h}} \frac{(\div \d, p_{T})}{||\d||_{\mathbf{H}^1} ||p_T||_{L^2}} \geq \beta_0 > 0 \quad 
\end{equation}
with $\beta_0$ independent of $h$. Here, we choose a Taylor-Hood type discretization with continuous, piecewise polynomials of degree $s \geq 2$ to approximate the displacement $\d$ over the whole domain $\Omega$ and continuous, piecewise polynomials of degree $s-1$ for the extracellular and intracellular total and fluid pressures over their respective subdomains. Thus, total and fluid pressure are discontinuous across the membrane interface $\Gamma$.
Letting $\mathcal{T}$ denote a triangulation of $\Omega$ which conforms to the interface $\Gamma$,
we define the discrete finite element spaces for $k \in \{ i,e \}$ as
\begin{equation}\label{eq:FE_spaces}
\begin{aligned}
\mathbf{V}_h &= \{\d \in \mathbf{C}(\Omega) \cap V: \d|_K \in \mathbb{P}_{s}(K)^d, \forall K \in \mathcal{T} \}, \\
Q_{T,h}^k &= \{q \in C(\Omega_k) \cap Q_T^k: q|_K \in \mathbb{P}_{s - 1}(K), \forall K \in \mathcal{T}\cap\Omega_k \}, \\
Q_{F,h}^k &= \{q \in C(\Omega_k) \cap Q_F^k: q|_K \in \mathbb{P}_{s - 1}(K), \forall K \in \mathcal{T}\cap\Omega_k \}
\end{aligned}
\end{equation}
and set $Q_{F,h} := Q_{F,h}^i \times  Q_{F,h}^e$ and $Q_{T,h} := Q_{T,h}^i \times  Q_{T,h}^e$. 
We note that it is also possible to employ polynomial order of $s$ for the fluid pressure space $Q_{F,h}$, which leads to an overall convergence rate of $O(h^{s})$ in the natural norms of each variable \cite{oyarzua_locking-free_2016}. However, due to the lower number of degrees of freedom and easier construction of multigrid methods for the coupled pressure block (see \Cref{sec:multigrid}), we use equal order spaces for total and fluid pressure.

The discrete stability of system \eqref{eq:total_pressure:discr} follows analogously to the proof of
its continuous counterpart in theorem \Cref{th:total_pressure_stability}, exploiting the inf-sup stability of the Taylor-Hood spaces.
We note that other Stokes-stable discretizations are equally applicable, see e.g. \cite{fortin1991mixed}.

\subsection{Efficient preconditioning of the pressure block involving $P_0^{-1}$}

The $Q$-block of the preconditioner $\mathcal{B}_d$ in \eqref{eq:TP_preconditioner_L20proj} includes the projection operator $P_0$, $P_0 q_T = q_T - \frac{1}{|\Omega|} \left(\int_\Omega q_T \,dx\right) \mathbf{1}$, which maps functions in $Q_{T,h}$ to their zero-mean component. We avoid the resulting dense block by expressing the projection as a low-rank correction as follows. Let $S_Q$ denote the full $2 \times 2$ matrix representation of this pressure block acting on coefficient vectors in $Q_{T,h} \times Q_{F,h}$.

The operator $P_0$ can be written as $\I - \Pi_0$, where $\Pi_0$ is the $L^2$-projection onto constant functions. When forming the matrix entries for the block corresponding to $Q_{T,h}$ (the $(1,1)$-block of $S_Q$), the contribution from $P_0 p_T$ becomes $(M_{TP} - \frac{1}{|\Omega|} \mathbf{m}\mathbf{m}^T)$. Here, $M_{TP}$ is the mass matrix for the $Q_{T,h}$ space, and $\mathbf{m}$ is the vector with entries $m_i = \int_\Omega \phi_i \,dx$ for each basis function $\phi_i \in Q_{T,h}$. Then, with the term $\frac{1}{|\Omega|} \mathbf{m}\mathbf{m}^T$ representing the matrix form of $\Pi_0$, the $2 \times 2$ pressure block matrix $S_Q$ from \eqref{eq:TP_preconditioner_L20proj} can be written as
\begin{equation}
S_Q = \begin{pmatrix} 
(\lambda^{-1}+1)M_{TP} - \frac{1}{|\Omega|}\mathbf{m}\mathbf{m}^T & K_{12} \\ 
K_{21} & K_{22} 
\end{pmatrix},
\end{equation}
where $K_{12}$ and $K_{21}$ are the matrix representations of the $-\alpha\lambda^{-1}I$ coupling terms, and $K_{22}$ is the matrix for the $((\alpha^{2}\lambda^{-1}+c_{0})I-\text{div }\kappa\nabla+L_{p}T)$ block acting on $Q_{F,h}$.

This structure allows us to express $S_Q$ as a rank-1 modification of a matrix $A$,
\begin{equation}
A := \begin{pmatrix} 
(\lambda^{-1}+1)M_{TP} & K_{12} \\ 
K_{21} & K_{22} 
\end{pmatrix}.
\end{equation}
Then, with the perturbation vector $\mathbf{y} = (\frac{1}{\sqrt{|\Omega|}}\mathbf{m}, \mathbf{0})^T$, where $\mathbf{0}$ is the zero vector in the $Q_{F,h}$ coefficient space, the pressure block becomes
\begin{equation}
S_Q = A - \mathbf{y}\mathbf{y}^T
\end{equation}
allowing us to avoid the dense (1,1) block of $S_Q$.
The action of $S_Q^{-1}$ required by the preconditioner can then be efficiently computed using the Sherman-Morrison-Woodbury (SMW) formula \cite{sherman1950adjustment}:
\begin{equation}\label{eq:smw}
(A-\mathbf{y}\mathbf{y}^T)^{-1} = A^{-1} + \frac{A^{-1}\mathbf{y}\mathbf{y}^T A^{-1}}{1-\mathbf{y}^T A^{-1}\mathbf{y}}.
\end{equation}
As the term $A^{-1}\mathbf{y}$ can be precomputed once during the setup-phase, the application of the preconditioner requires only a single evaluation of $A^{-1}$ (and a few vector products), making the preconditioner \eqref{eq:TP_preconditioner_L20proj} only slightly more expensive than its $P_0$-free counterpart \eqref{eq:TP_preconditioner}. Finally, we note that the action of $A^{-1}$ in the SMW formula can either be computed exactly with Cholesky decomposition, or approximated with AMG.

\subsection{Multigrid approximations}\label{sec:multigrid}

To construct an efficient and scalable block diagonal preconditioner of the form \eqref{eq:TP_preconditioner}, we need to replace the exact inverses in each block with computationally inexpensive, yet spectrally equivalent operators. The displacement block is a standard second-order elliptic operator, and we apply a single AMG V-cycle with Hypre's BoomerAMG solver \cite{yang2002boomeramg,hypre} with default settings, except for a larger \textit{strong threshold} ($\theta=0.5$) and an additional smoother applications per V-cycle (three instead of default two).

\subsubsection{Parameter robust AMG for the $Q$-block}

The $2 \times 2$-block of total and fluid pressure in both preconditioner \eqref{eq:TP_preconditioner} and \eqref{eq:TP_preconditioner_L20proj} is not standard, and can potentially pose a challenge for AMG-based parameter-robust approximations. Considering the case of mixed boundary conditions on the displacement, we decompose the inner-product operator on $Q_{T,h} \times Q_{F,h}$ as follows
\begin{align}
\begin{bmatrix}
    I & 0 \\
    0 &  c_0I - \div \kappa \nabla
\end{bmatrix} + \frac{1}{\lambda}
\begin{bmatrix}
    I & - \alpha I \\
    - \alpha I & \alpha^2 I 
\end{bmatrix}+ L_p
\begin{bmatrix}
    0 & 0 \\
    0 & T 
\end{bmatrix}.
\end{align}\label{eq:q_block}
Here, the second term represents the coupling between the total pressure $p_T$ and the fluid pressure $p_F$, and can become singular in certain parameter regimes. For instance, setting $\alpha=1$, we find the vector $(1, 1)^T$ in the kernel of the coupling operator, and note that the strength of the singular perturbation is controlled by $\lambda^{-1}$.
The third term represents the interfacial coupling between the intra- and extracellular fluid pressure, controlled by the membrane permeability $L_p$. Similar to the bulk coupling term, interfacial coupling term features a non-trivial kernel, namely any (0, $p_F$) for $p_F \in Q_F$ with $(\jump{p_F}, \jump{p_F})_{\Gamma}=0$.

Due to the presence of these singular perturbations, it is not clear a apriori whether a standard AMG approach yields sufficiently robust approximations. 
Specific multilevel methods for nearly-singular systems are available in the literature (see e.g. \cite{lee2007robust} for bulk coupling and \cite{budivsa2024algebraic} for interface perturbations), and suggest that careful treatment of the (near) kernel is crucial for parameter-robust AMG methods for singularly perturbed elliptic operators \cite{budivsa2024algebraic, lee2007robust}. However, such methods are only partly directly available in PETSc and thus increase implementation difficulty,   often hindering their adoption in application code.
Here, we speculate that a standard unknown-based AMG approach might still deliver sufficiently robust performance, and test this hypothesis for our specific case. To this end we employ the conjugate gradient (CG) method preconditioned with a single AMG V-cycle to the inner product operator of $Q_{T,h} \times Q_{F,h}$, and compare different AMG approximations with respect to the spectral condition number and the number of CG iterations required.

Motivated by the analysis in \cite{budivsa2024algebraic} and \cite{lee2007robust}, we include BoomerAMG's nodal
solution strategy~\footnote{We activate both nodal coarsening and nodal relaxation via the PETSc flags \textit{pc\_hypre\_boomeramg\_nodal\_coarsen} and \textit{pc\_hypre\_boomeramg\_nodal\_relaxation}.} in addition to its default unknown-based approach in our comparison. Further, we identify the \textit{strong threshold} $\theta \in [0,1]$ as a key parameter for AMG performance, and test  $\theta \in \{0.3, 0.5, 0.7 \}$.
In \Cref{tb:pressure_coupling} we observe that the nodal approach indeed performs substantially better in the strong coupling regime ($\lambda \ll 1$), but is negatively affected by an increase in $\kappa$. In contrast, the default unknown-based approach provides good approximations as long as $\lambda \geq 1$, with a small advantage for lower strong threshold values.
While these results generally confirm the crucial role of kernel-capturing smoothers, we note that $\lambda \geq 1$ holds for the broad class of materials with Poisson's ratio between 0 and 0.5, rendering the singular perturbation from the fluid and total pressure coupling less problematic.

Next, we investigate the role of the interfacial perturbation in AMG approximations. \Cref{tb:Lp_stability} shows
that the robustness with respect to $L_p$ highly depends on the strong threshold. Full robustness is only achieved
for $\theta=0.7$, while lower threshold values deteriorate in the $L_p \gg \kappa$ regime, i.e. when singular perturbation
is dominant. Thus, we proceed with the strong threshold value $\theta=0.7$ for the remainder of the paper. 

\begin{table}
\setlength{\tabcolsep}{4pt}
\scriptsize
\centering
\begin{tabular}{l|ccc|ccc}
\toprule
 & \multicolumn{3}{c}{unknown-based (default)} & \multicolumn{3}{c}{nodal} \\
 \midrule
\diagbox{$(\lambda, \kappa)$}{$\theta$} & 0.30 & 0.50 & 0.70 & 0.30 & 0.50 & 0.70 \\
\midrule
($10^{-5}$, $10^{-7}$) & 77869.50 (85) & 77903.03 (85) & 77894.26 (85) & 1.90 (12) & 1.90 (13) & 1.90 (13) \\
($10^{-5}$, $10^{-3}$) & 21450.65 (520) & 21267.23 (514) & 21549.41 (525) & 5.97 (23) & 6.21 (24) & 6.20 (23) \\
($10^{-5}$, $10^{0}$) & 7225.76 (483) & 6892.63 (479) & 7246.39 (481) & 4451.67 (236) & 4522.61 (234) & 4228.18 (233) \\
($10^{-5}$, $10^{3}$) & 87008.32 (77) & 85399.73 (78) & 93670.99 (98) & 62450.49 (78) & 70762.48 (85) & 71540.24 (81) \\
\midrule
($10^{0}$, $10^{-7}$) & 2.32 (14) & 2.24 (14) & 2.28 (14) & 1.46 (10) & 1.36 (9) & 1.30 (9) \\
($10^{0}$, $10^{-3}$) & 1.64 (10) & 1.65 (10) & 1.65 (10) & 2.70 (12) & 2.68 (12) & 1.69 (11) \\
($10^{0}$, $10^{0}$) & 2.03 (12) & 2.11 (13) & 2.36 (14) & 14.85 (22) & 2.06 (14) & 4.20 (17) \\
($10^{0}$, $10^{3}$) & 2.00 (13) & 2.01 (13) & 2.73 (15) & 12404.55 (34) & 22.44 (17) & 4767.14 (27) \\
\midrule
($10^{5}$, $10^{-7}$) & 1.12 (7) & 1.27 (8) & 1.35 (9) & 10.25 (19) & 10.16 (22) & 4.08 (18) \\
($10^{5}$, $10^{-3}$) & 1.59 (9) & 1.68 (10) & 2.72 (14) & 583.62 (26) & 27.14 (17) & 164.14 (22) \\
($10^{5}$, $10^{0}$) & 1.65 (9) & 1.69 (11) & 2.74 (14) & 563742.81 (37) & 915.96 (17) & 216598.72 (30) \\
($10^{5}$, $10^{3}$) & 1.65 (9) & 1.68 (11) & 2.74 (14) & $1.44 \cdot 10^8$ (37) & 914670.74 (20) & 
$7.79 \cdot 10^7$ (26) \\
\bottomrule
\end{tabular}
\label{tb:pressure_coupling}
\caption{Spectral condition number of the Riesz map with respect to the inner product  on $Q_{T,h} \times Q_{F,h}$ induced by \eqref{eq:TP_preconditioner} when using different AMG methods for preconditioning. The number in parenthesis denotes the number of CG iterations needed to reach a relative residual of $10^{-10}$. The computations are carried out on a astrocyte geometry with 2.5\,\textmu m sidelength (see \Cref{fig:astro_meshes}), and the remaining parameters are $c_0=10^{-6}$, $\alpha=1$ and $L_p=0$.}
\end{table}

\begin{table}
\setlength{\tabcolsep}{5pt}
\centering
\scriptsize
\begin{tabular}{l|ccc}
\toprule
 & \multicolumn{3}{c}{unknown-based (default)} \\
  \midrule
\diagbox{$(L_p, \kappa)$}{$\theta$} & 0.30 & 0.50 & 0.70 \\
\midrule
($10^{-9}$, $10^{-7}$) & 2.31 (24) & 2.33 (24) & 2.33 (24) \\
($10^{-9}$, $10^{-3}$) & 1.64 (17) & 1.65 (17) & 1.65 (17) \\
($10^{-9}$, $10^{0}$) & 2.05 (22) & 2.11 (22) & 2.36 (25) \\
($10^{-9}$, $10^{3}$) & 2.01 (21) & 2.01 (23) & 2.73 (27) \\
\midrule
($10^{-2}$, $10^{-7}$) & 1.68 (12) & 1.64 (12) & 1.63 (12) \\
($10^{-2}$, $10^{-3}$) & 1.24 (9) & 1.25 (9) & 1.25 (9) \\
($10^{-2}$, $10^{0}$) & 2.05 (13) & 2.09 (13) & 2.23 (15) \\
($10^{-2}$, $10^{3}$) & 2.01 (14) & 2.01 (14) & 2.91 (16) \\
\midrule
($10^{2}$, $10^{-7}$) & 5019.55 (502) & 5219.82 (163) & 3.54 (19) \\
($10^{2}$, $10^{-3}$) & 41.49 (69) & 21.80 (49) & 1.71 (12) \\
($10^{2}$, $10^{0}$) & 2.06 (9) & 2.12 (10) & 2.35 (12) \\
($10^{2}$, $10^{3}$) & 2.01 (9) & 2.01 (11) & 1.92 (15) \\
\bottomrule
\end{tabular}
\label{tb:Lp_stability}
\caption{Spectral condition number of the Riesz map with respect to the inner product  on $Q_{T,h} \times Q_{F,h}$ induced by \eqref{eq:TP_preconditioner} when using different AMG strong thresholds for preconditioning. The number in parenthesis denotes the number of CG iterations needed to reach a relative residual of $10^{-10}$. The computations are carried out on a 2.5\,\textmu m sidelength astrocyte geometry (see \Cref{fig:astro_meshes}), and the remaining parameters are $c_0=10^{-6}$, $\alpha=1$ and $\lambda=1$.}
\end{table}

\begin{remark}[Spatially varying material parameters]
\label{rem:heterogeneity}
For simplicity, we treat the material parameters as scalar constants in the numerical experiments to reduce the dimension of the parameter space. However, with the exception of the shear modulus $\mu$ (due to the rescaling in \eqref{eq:total_pressure}), the theoretical framework applies equally to spatially varying parameters. The weighted norms naturally incorporate spatial inhomogeneity. Consequently, the condition number of the preconditioned system remains bounded independently of the magnitude of parameter jumps. 

In terms for practical realization we note that AMG methods can handle spatial heterogeneity \cite{xu2017algebraic}. 
However, in case their performance deteriorates (for example, due to presence of large jumps) other approximations of the 
preconditioner blocks can be used, e.g. domain decomposition methods such as GenEO \cite{spillane2014abstract} or BDDC \cite{badia2016physics}.
\end{remark}

\section{Discrete preconditioners for the cell-by-cell model}\label{sec:numerical_results}

In this section, we present numerical results for the preconditioner and solver introduced in the previous section. We demonstrate parameter-robustness of the proposed preconditioners using exact inverses of each diagonal block in 2D and 3D.
In addition, we investigate the performance of the suggested multigrid methods and showcase their efficiency
on a large-scale example of cellular swelling in the rat visual cortex.
The numerical experiments are carried out the finite element software FEniCS~\cite{alnaes2015fenics} with the extension multiphenics~\cite{ballarin2019multiphenics} to handle coupling of function spaces defined on subdomains.

\subsection{Approximation error}

We begin by verifying the approximation properties of the scheme by choosing the right-hand-side of system \cref{eq:total_pressure} such that the exact solution on the unit square $\Omega=(0,1)^2$ with the interface $\Gamma$ at $x=0.5$ is given by
\begin{equation}
\begin{aligned}
\d(x, y) &= (2 \pi \sin{(2 \pi x)}  \cos{(2 \pi y)}, -2 \pi \cos{(2 \pi x)}  \sin{(2 \pi y)} )^T, \\
p_i(x, y) &= \sin{(\pi x)}   \cos{(3.4 \pi y)}, \\
p_e(x, y) &= 1 + p_i(x, y).
\end{aligned}
\end{equation}
Using discretization by \eqref{eq:FE_spaces} with $s=2$ we observe quadratic convergence of
$\d_h$, $p_{T, h}$ in their natural norms ($H^1$, respectively $L^2$) while the $H^1$-norm of
$p_F - p_{F, h}$ decreases linearly with $h$. We recall that $Q_{F, h}$ is discretized by linear elements.
The observed rates are inline with theoretical results of \cite{oyarzua_locking-free_2016}
which bound the convergence rates in terms of the minimal polynomial order of the pressure
and displacement spaces (i.e. order 1 in our case).

\begin{table}[]
  \centering
\setlength{\tabcolsep}{4pt}
\small
\begin{tabular}{lccc}
\toprule
$N$ & $||\d - \d_h||_{H^1}$ & $||p_F - p_{F, h}||_{H^1}$  & $|| p_T -  p_{T, h}||_{L^2}$ \\
\midrule
8   & 0.82849 (-)    & 2.42620 (-)    & 0.05144 (-)                  \\
16  & 0.20997 (1.98) & 1.22622 (0.98) & 0.01260 (2.03)  \\
32  & 0.05273 (1.99) & 0.61482 (1.00) & 0.00313 (2.01)  \\
64  & 0.01320 (2.00) & 0.30763 (1.00) & 0.00078 (2.00)  \\
128 & 0.00330 (2.00) & 0.15384 (1.00) & 0.00020 (2.00) \\
\bottomrule
\end{tabular}
\label{tab:EOC}
\caption{Approximation errors of the Taylor-Hood type 
finite element discretization \eqref{eq:FE_spaces} with $s=2$ on the domain $\Omega=(0,1)^2$ with $\Gamma$ at $x=0.5$, uniformly triangulated in $2 N^2$ cells. The estimated order of convergence is shown in the parenthesis.}
\end{table}

\subsection{Parameter ranges for sensitivity analysis}\label{sec:parameter_values}
Next, we identify the parameter regimes for which we will assess robustness of the preconditioners.
We chose the parameter ranges based on a scaling analysis and several real-world applications. To this
end let us non-dimensionalize the static version of problem \eqref{eq:full_biot} arising from time
discretization with time step size $\tau > 0$ by introducing the non-dimensional variables 
$\tilde{\d} = \d / d_0$,
$\tilde{p} = p/ p_0$,
and the differential operators 
$\tilde{\div} = L \div$,
$\tilde{\nabla} = L \nabla$,and
$\tilde{\eps} = L \eps$, where 
$d_0$, $p_0$, and $L$ are characteristic displacement, pressure, and length scale, respectively.
Inserting the non-dimensional variables yields a static problem of the type
\begin{align}
\begin{aligned}
- L^{-1} \tilde{\div} [ 2 \mu L^{-1} d_0\tilde{\eps} (\tilde{\d}) + \lambda L^{-1} d_0 \tilde{\div}  \tilde{\d} I  - \alpha p_0 \tilde{p} I] &= 0, \\
c_0 \tau^{-1} p_0 \tilde{p} + \alpha L^{-1} d_0 \tau^{-1} \tilde{\div} \tilde{\d} - L^{-2} p_0 \tilde{\div} (\kappa \tilde{\nabla} \tilde{p})  &= 
\tilde{g},
 \end{aligned}
\end{align}
where $\tilde{g}$ contains additional contributions from time discretization.
Multiplying the momentum equation with $\frac{L^{2}}{2 \mu d_0}$ and the mass conservation equation with $\frac{p_0 \tau L^2}{2 \mu d_0^2}$, we find
\begin{align}
\begin{aligned}
-\tilde{\div} [\tilde{\eps}(\tilde{\d}) +
 \frac{\lambda}{2 \mu} \tilde{\div}  \tilde{\d} I  - \frac{\alpha p_0 L}{2 \mu d_0} \tilde{p} I] &= 0, \\
\frac{c_0 p_0^2 L^2}{2 \mu d_0^2} \tilde{p}
+ \frac{\alpha p_0 L}{2 \mu d_0}  \tilde{\div} \tilde{\d} 
- \frac{p_0^2 \tau}{2 \mu d_0^2} \tilde{\div} (\kappa \tilde{\nabla} \tilde{p})  &= \frac{p_0 \tau L^2}{2 \mu d_0^2} \tilde{g}
 \end{aligned}
\end{align}
with the non-dimensionalized interface coupling condition
\begin{equation}
    \frac{L_p L}{\kappa} \jump{\tilde{p}} = - \tilde{\nabla} \tilde{p} \cdot \n.
\end{equation}

This allows us to identify the following dimensionless numbers, which effectively characterize our system: The dimensionless Darcy number $Da$, the dimensionless storage coefficient $S$, the dimensionless Biot-Willis coefficient $BW$, the dimensionless elastic modulus $E$, and the dimensionless coupling coefficient $Cp$, defined as 
\begin{equation}
Da = \frac{\kappa p_0^2 \tau}{2 \mu d_0^2},\quad
S = \frac{c_0 p_0^2 L^2}{2 \mu d_0^2} ,\quad
BW = \frac{\alpha p_0 L}{2 \mu d_0},\quad
E = \frac{\lambda}{2 \mu},\quad
Cp = \frac{L L_p}{\kappa},\quad
\end{equation}
such that
\begin{align}
\begin{aligned}
-\tilde{\div} [\tilde{\eps}(\tilde{\d}) +
 E \tilde{\div}  \tilde{\d} I  - BW \tilde{p} I] &= 0 \quad &\text{ in }\tilde{\Omega},\\
S \tilde{p}
+ BW \tilde{\div} \tilde{\d} 
- Da \tilde{\div} ( \tilde{\nabla} \tilde{p})  &= 0 &\text{ in }\tilde{\Omega},\\
Cp \jump{\tilde{p}} + \tilde{\nabla} \tilde{p} \cdot \n &= 0 &\text{ on }\tilde{\Gamma}.
 \end{aligned}
\end{align}
Comparing the non-dimensional system with the simplified system \eqref{eq:strong}, we find that we can effectively interpret $Da$ as $\tilde{\kappa}$, $S$ as $\tilde{c_0}$, $BW$ as $\tilde{\alpha}$, $E$ as $\tilde{\lambda}$, and $Cp \cdot Da$ as $\tilde{L_p}$ for unit scaled systems.

We remark that this choice of scaling is general, and assumes that the characteristic scales of displacement, pressure, length and time are independent. However, in many practical applications, further simplifications are possible. For instance, in the absence of a strong forcing term or dominating boundary displacements, the solid stress terms are balanced by the fluid pressure, implying that the characteristic displacement obeys $d_0 = \frac{\alpha p_0 L}{2 \mu}$, and thus $BW=1$, $S=2 \mu c_0 / \alpha^2$ and $Da=\frac{2 \mu \tau \kappa}{\alpha^2 L}$. 

Next, we consider three real-world examples to determine the relevant parameter ranges.
%
%\subsection{Relevant parameter ranges in applications}\label{sec:parameter_values}
%
Our work is primarily motivated by the problems arising from the biomechanics of osmotically induced cellular swelling in the mammalian brain. However, similar problems also arise in various other fields, where a thin, less permeable layer separates two or more poroelastic materials. Here, we provide three examples of potential real-world applications across different scales and estimate relevant parameter ranges.

\begin{example}[Cellular swelling in the mammalian brain]\label{ex:1}
In the brain, the composition of interstitial fluid varies considerably between physiological states such as sleep and locomotion and is disturbed in pathological states \cite{rasmussen2020interstitial}. These changes induce substantial osmotic forces on the cells, raising important questions about volume regulation and cellular fluid exchange. The cells can be described as poroelastic materials, separated from the poroelastic extracellular matrix via their cell membrane.
In this scenario, we assume a characteristic length of $20\,\mu$m, a characteristic displacement of $100\,$nm, a characteristic pressure of $10\,$Pa, and a characteristic time step size of $0.1\,$s (see e.g. \cite{causemann2025stretch}). Together with typical material parameters specified in \Cref{tb:parameter_values}, we find the ranges $Da \in [10^{-2}, 10^2]$, $S \in [10^{-5}, 1]$, $BW \in [1, 10]$, $E \in [10, 10^4]$, $Cp \in [10^{-8}, 10^{-2}]$ and $Da \cdot Cp \in [10^{-7}, 10^{-3}]$ for the non-dimensional groups. \\
\end{example}

\begin{example}[Tissue engineering]\label{ex:2}
In a tissue engineering, scenario, a hydrogel scaffold (a poroelastic medium) is split into multiple parts by a semi-permeable membrane, mimicking biological tissues with limited fluid flow and cell migration across a membrane (see e.g. \cite{ozkendir2024engineering}). Concretely, we consider two alginate hydrogels ($E=0.05\,$MPa, $\nu=0.38$, $\kappa = 7 \cdot 10^{-18}\,\text{m}^2$), separated by a semipermeable membrane ($L_p \in [ 10^{-12}, 10^{-16}]$), and are interested in flow and displacement over a length scale of $5\,$mm, computed with a time step size of 1 hour. Here, we obtain the ranges $Da \in [10^{-7}, 10^{-6}]$, $S \in [10^{-6}, 10^{-3}]$, $BW \in [10^{-2}, 10^{-1}]$, $E \in [10, 10^2]$, $Cp \in [10^{-4}, 10]$ and $Da \cdot Cp \in [10^{-10}, 10^{-5}]$.\\
\end{example}

\begin{example}[Aquifer systems]\label{ex:3}
In hydrology, aquifer systems are commonly modeled as poroelastic materials. For instance, an unconfined aquifer (a poroelastic medium) can be separated from a confined aquifer (another poroelastic medium) by an aquitard or aquiclude (a thin layer of less permeable rock).
Concretely, we consider two ground-water saturated layers of sand ($k=10^{-9}\,\text{m}^2$), separated by a thin layer of clay ($L_p \in [10^{-18}, 10^{-14}] \text{m}\,\text{Pa}^{-1}\text{s}^{-1}$), represented as the two-dimensional interface between the two sand layers.
Assuming relevant length and time step scales of $L=500\,$m and $\tau=1\,$day, and considering typical material parameters (see \Cref{tb:parameter_values}), results in the dimensionless groups  $Da \in [10^{-7}, 10^3]$, $S \in [10^{-3}, 10]$, $BW \in [10^{-2}, 10]$, $E \in [10, 10^5]$, $Cp \in [10^{-9}, 10^{2}]$ and $Da \cdot Cp \in [10^{-9}, 10^{-2}]$.\\
\end{example}

\begin{table}[!ht]
\centering
\small
\begin{tabular}{p{1.4cm}p{2.1cm}p{1.4cm}p{1.6cm}p{1.45cm}p{1.5cm}p{2.5cm}}
\toprule
\textbf{Scenario} & \textbf{Intrinsic permeability ($K$)} & \textbf{Biot-Willis coefficient ($\alpha$)} & \textbf{Storage coefficient ($c_0$)} & \textbf{Young modulus ($E$)} & 
\textbf{Poisson's ratio ($\nu$)} & \textbf{Membrane conductivity ($L_p$)} \\ \midrule
\textbf{Cellular swelling}  & $10^{-16} - 10^{-14} \, \text{m}^2$ & 1.0 & $10^{-8} - 10^{-5} \, \text{Pa}^{-1}$ & $500 - 1500 \, \text{Pa}$& $0.17 - 0.48$ & $10^{-14} - 10^{-11} \, \text{m s}^{-1}\text{Pa}^{-1}$ \\
\midrule
\textbf{Tissue engineering}   & $10^{-18} \, \text{m}^2$ & 1.0 & $10^{-6} - 10^{-4} \, \text{Pa}^{-1}$ & $5 \cdot 10^4 \, \text{Pa}$ & $0.38$ & $10^{-16} - 10^{-12} \, \text{m s}^{-1}\text{Pa}^{-1}$ \\
\midrule
\textbf{Aquifer Systems}  & $10^{-9} - 10^{-16} \, \text{m}^2$ & 0.6 - 1.0 & $10^{-9} - 10^{-11} \, \text{Pa}^{-1}$ & $10^{9} - 10^{10} \, \text{Pa}$ & $0.15 - 0.35 $& $10^{-16} - 10^{-14} \, \text{m s}^{-1}\text{Pa}^{-1}$ \\ \bottomrule
\end{tabular}
\label{tb:parameter_values}
\caption{Parameter ranges for various poroelastic scenarios. The hydraulic conductivity $\kappa$ is obtained from the intrinsic permeability $K$ via the relation $\kappa = K/\mu_f$, assuming a fluid viscosity of $\mu_f = 10^{-3}\,\text{Pa}\cdot\text{s}$ (water) in all scenarios. The Lamé parameters $\lambda$ and $\mu$ are computed from Young's modulus $E$ and Poisson's ratio $\nu$ via $\mu=E/(2(1+\nu))$ and $\lambda=E\nu/((1+\nu)(1-2\nu))$. The estimated ranges for each scenario are based on the following references: Cellular Swelling \cite{jin2016spatial,mokbel2020poisson,lu2006viscoelastic,causemann2025stretch}; Tissue engineering: \cite{hollister2005porous,ozkendir2024engineering,o2011biomaterials};
Aquifer Systems: \cite{zimmerman2000micromechanics}.
}
\end{table}

\begin{table}[h!]
\centering
\small
\begin{tabular}{p{1.4cm}p{1.05cm}p{0.75cm}p{1.4cm}p{1.6cm}p{0.9cm}p{0.9cm}p{0.9cm}p{0.9cm}p{0.9cm}}
\toprule
\textbf{Scenario} & \textbf{length Scale (L)} & \textbf{time step ($\tau$)} & \textbf{char. pressure ($p_0$)} & \textbf{char. displacement ($d_0$)} & $\mathbf{Da}$ & 
\textbf{$S$} & \textbf{$BW$} & \textbf{$E$} & \textbf{$Cp$} \\ \midrule
\textbf{Cellular swelling} & $20\,\mu$m & $0.1\,$s &  $10\,$Pa & $100\,$nm & $10^{-2}$--$10^2 $ & $10^{-5}$--$ 1$ & $1$--$ 10$ & $10$--$ 10^4$ & $10^{-8}$--$ 10^{-2}$  \\
\midrule
\textbf{Tissue engineering} & $5\,$mm & $1\,$h &  $10\,$Pa & $0.1\,$mm & $10^{-7}$--$ 10^{-6}$ 
& $10^{-6}$--$10^{-3}$ & $10^{-2}$--$ 10^{-1}$ & $10$--$ 10^2$ & $10^{-4}$--$ 10$  \\
\midrule
\textbf{Aquifer Systems} & $500\,$m & $24\,$h & $1\,$MPa & $1\,$m & $10^{-7} $--$ 10^3$ & $10^{-3}$--$ 10$ & $10^{-2}$--$ 10$ & $10$--$ 10^5$ & $10^{-9}$--$ 10^{2}$ \\
\bottomrule
\end{tabular}
\caption{Typical ranges for the dimensionless groups $Da$, $S$, $BW$, $E$, and $Cp$ in
  real-world applications \Cref{ex:1,ex:2,ex:3}.}
\label{tb:dim_groups}
\end{table}

In light of the obtained ranges for the dimensionless groups in real-world applications, and noting that they can effectively represent the parameters in system \eqref{eq:full_biot}, we deem the ranges $\alpha \in [10^{-2}, 1^{2}]$, $\kappa \in [10^{-7}, 10^{3}]$,  $\lambda \in [10, 10^{5}]$, $c_0 \in [10^{-6}, 10]$ and $L_p \in [10^{-9}, 10^{-2}]$ suitable to cover a broad range of potential applications. Restricting the parameter space for our numerical experiments, we focus on the most challenging case with $c_0=10^{-6}$.

\subsection{Two-dimensional examples}

\begin{figure}
    \centering
\includegraphics[trim={0.8cm 0 3cm 0},clip, height=0.48 \textwidth]{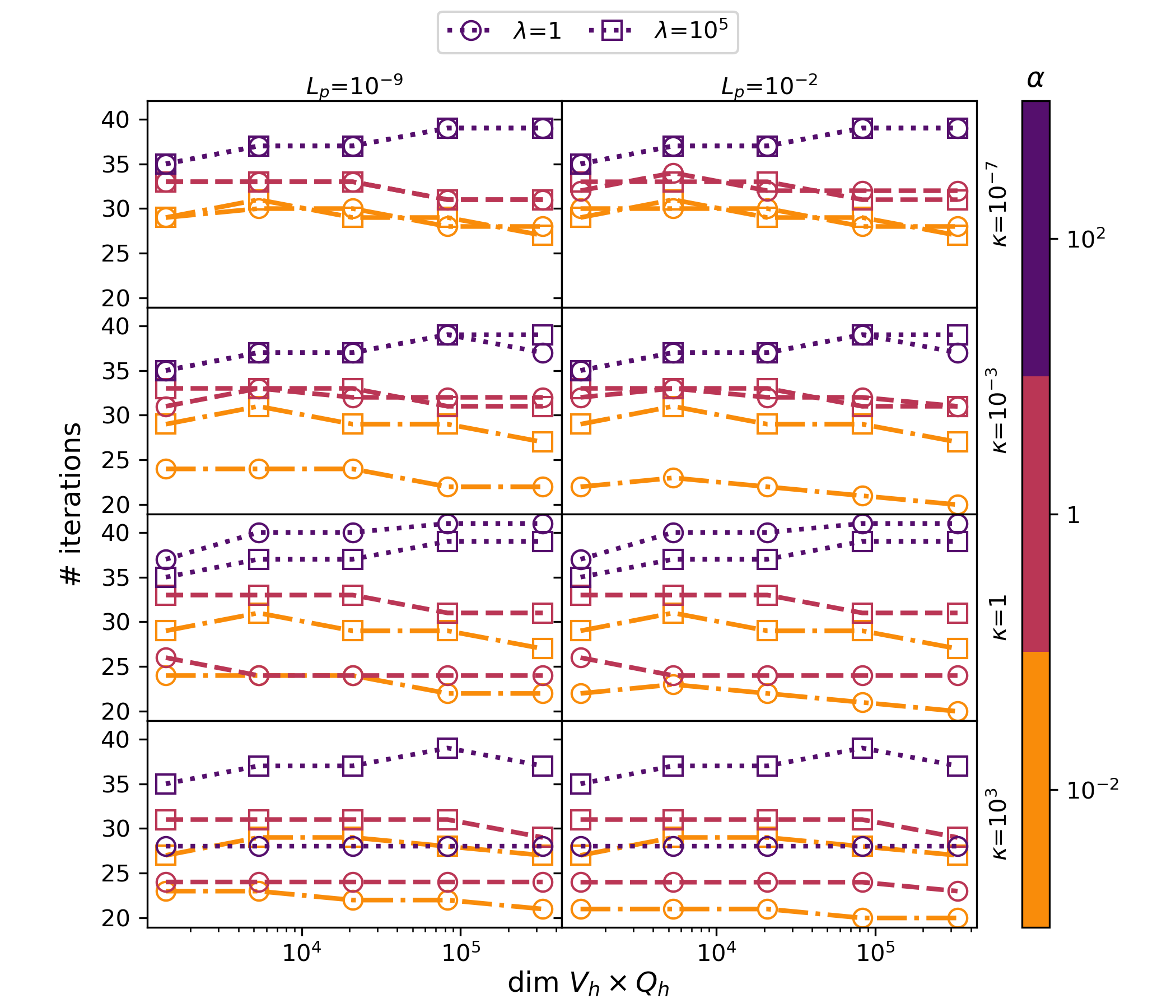}
\includegraphics[trim={0.2cm 0 0.5cm 0},clip, height=0.48 \textwidth]{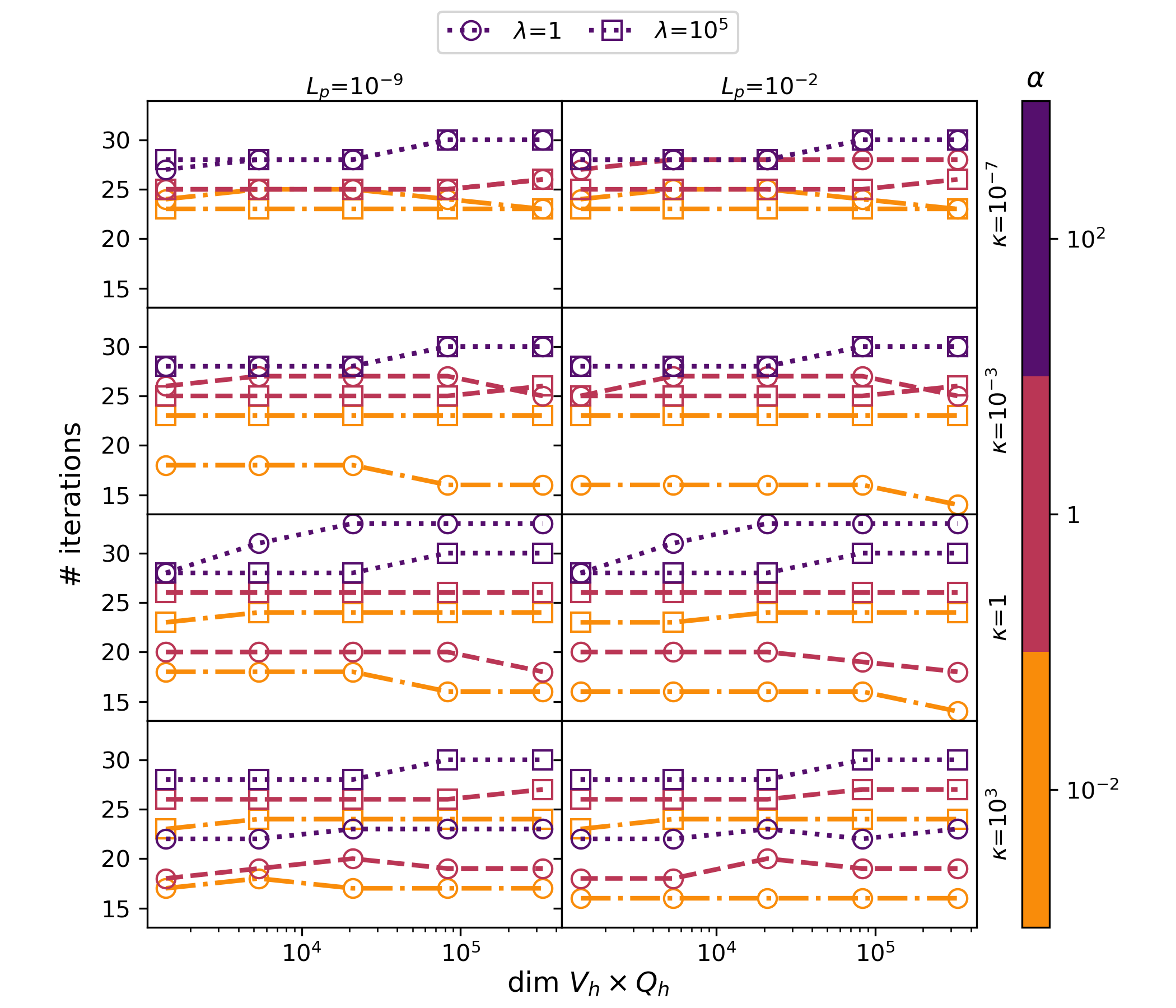}
\caption{Number of MinRes iterations for test cases on a unit square with exact inverses of the preconditioner blocks computed via Cholesky factorization: Mixed displacement boundary conditions ($|\Gamma_d| >0$ and $|\Gamma_t| >0$) with preconditioner \eqref{eq:TP_preconditioner} (left); and full displacement Dirichlet boundary conditions ($\partial \Omega = \Gamma_d$) with preconditioner \eqref{eq:TP_preconditioner_L20proj} (right).}
\label{fig:square}
\end{figure}

First, we demonstrate the parameter and mesh independent performance of the preconditioner \eqref{eq:TP_preconditioner} for mixed displacement boundary conditions ($|\Gamma_d| >0$ and $|\Gamma_t| >0$) and preconditioner \eqref{eq:TP_preconditioner_L20proj} for full Dirichlet displacement boundary conditions ($\partial \Omega = \Gamma_d$) by solving the preconditioned system with manufactured solutions on a unit square geometry, divided in an extracellular domain $\Omega_e = (0, 0.5) \times (0,1)$ and an intracellular domain  $\Omega_e = (0.5, 1) \times (0,1)$.
In both cases, we impose full Neumann pressure boundaries ($\Gamma_f=\partial \Omega$), and solve the system up to relative tolerance of $10^{-10}$. The individual preconditioner blocks are inverted exactly via Cholesky factorization with the MUMPS direct solver package \cite{MUMPS}.
Both preconditioners yield stable iteration counts across all parameter regimes and mesh sizes, varying between 14 and 39 (see \Cref{fig:square}).

\subsection{Three-dimensional examples and multigrid approximations}

Next, we test the performance of the multigrid preconditioner presented in \Cref{sec:multigrid} on a unit cube geometry (again divided into an extra- and intracellular domain at $x=0.5$) and compare the number of MinRes iterations with the exact counterparts computed via Cholesky factorization of the individual blocks. 

\begin{figure}
\centering
\includegraphics[trim={0.8cm 0 3cm 0},clip, height=0.48 \textwidth]{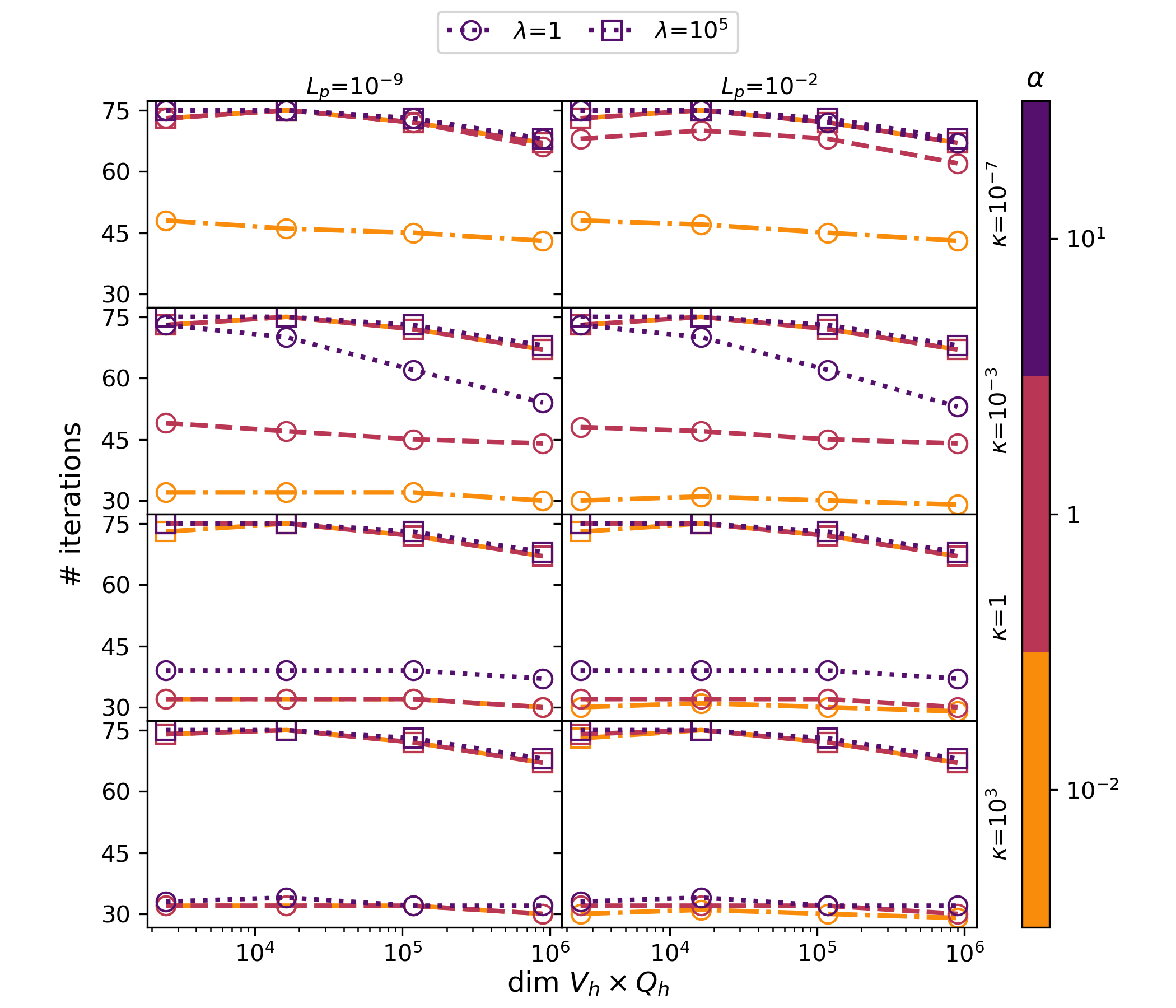}
\includegraphics[trim={0.2cm 0 0.5cm 0},clip, height=0.48 \textwidth]{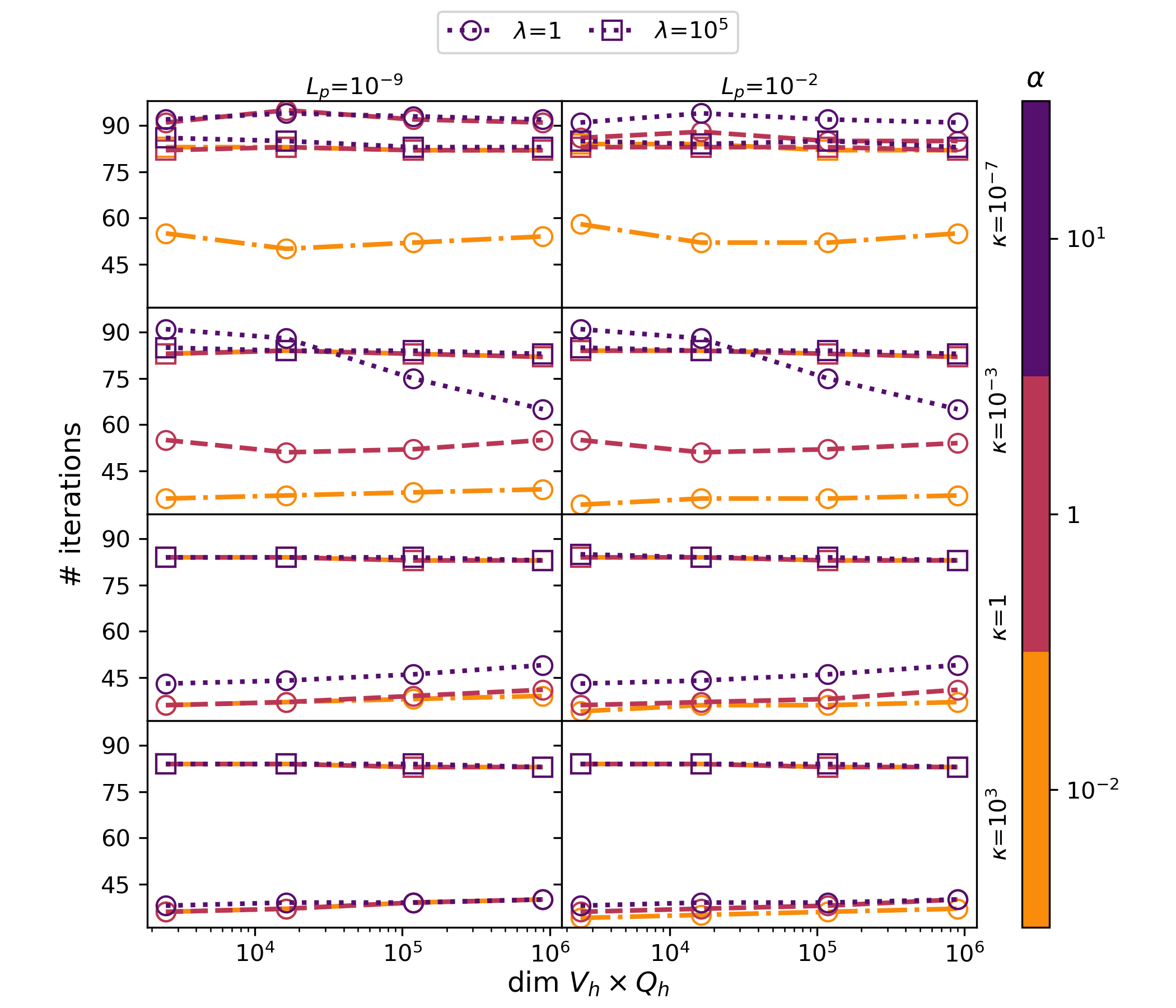}
\caption{Number of MinRes iterations for test cases on a unit cube with mixed displacement boundary conditions ($|\Gamma_d| >0$ and $|\Gamma_t| >0$) and preconditioner \eqref{eq:TP_preconditioner}, with exact inverses of the individual blocks (left), and algebraic multigrid approximations (right).}
\label{fig:cube_mixed}
\end{figure}

Starting with the mixed displacement boundary condition case, we find that the number of MinRes iterations remains stable under mesh refinement across all parameter regimes. The multigrid preconditioner leads to a small (maximum 20\%) increase compared to the exact preconditioner, ranging between 33 and 93 iterations until a relative residual of $10^{-10}$ is reached (see \Cref{fig:cube_mixed}).

\begin{figure}
\centering
\includegraphics[trim={0.8cm 0 3cm 0},clip, height=0.48 \textwidth]{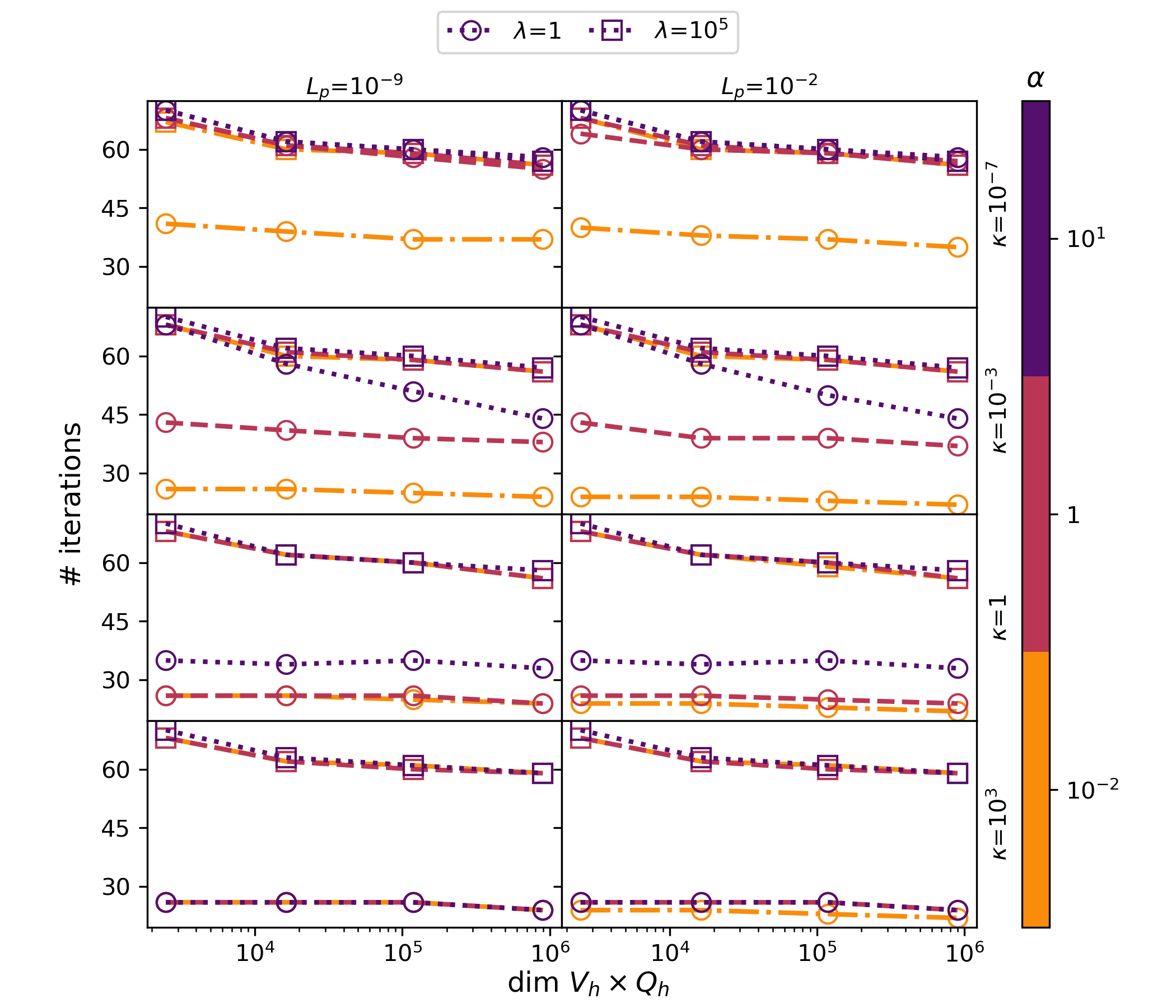}
\includegraphics[trim={0.2cm 0 0.5cm 0},clip, height=0.48 \textwidth]{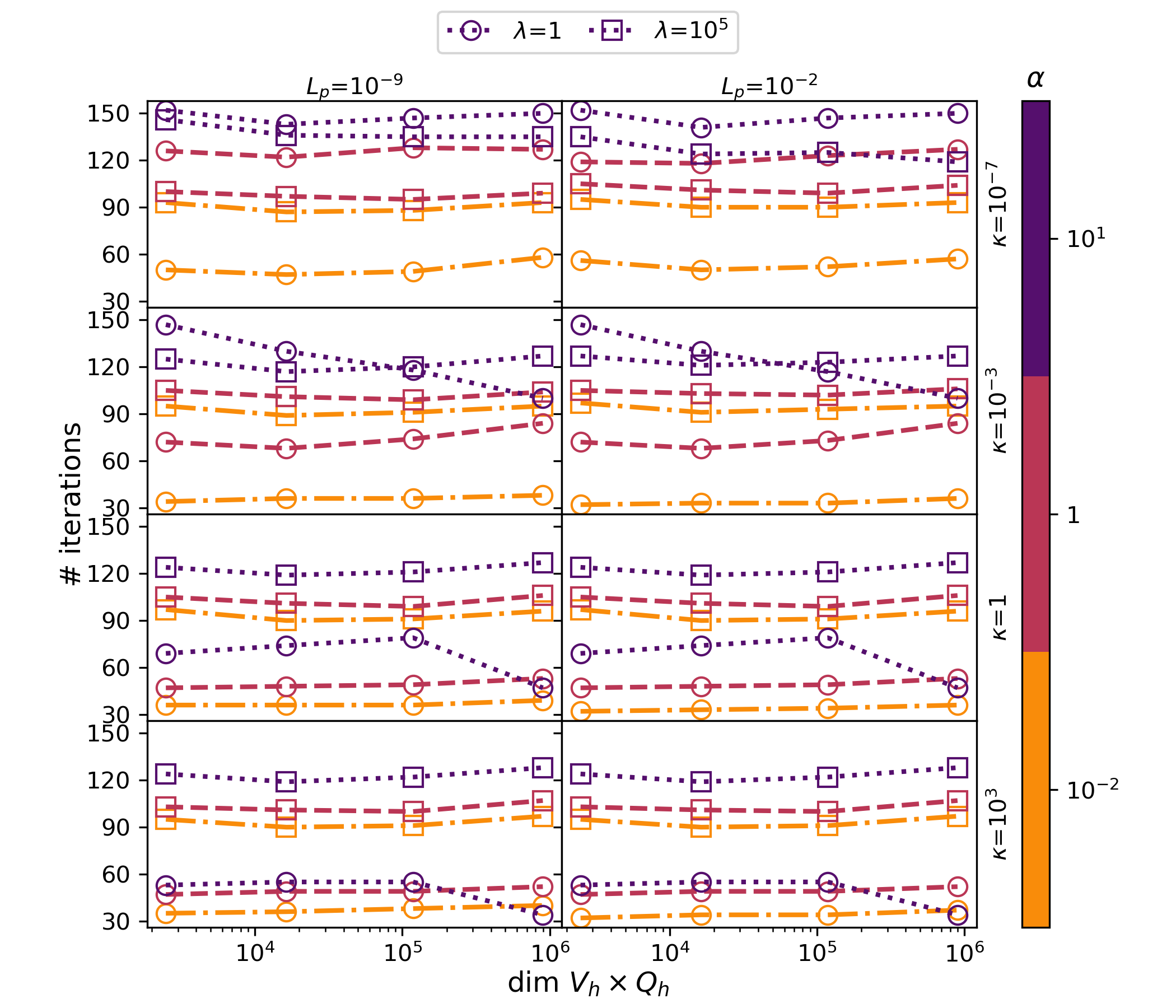}
\caption{Number of MinRes iterations for test cases on a unit cube with full Dirichlet displacement boundary conditions ($\Gamma_d = \partial \Omega$) and preconditioner \eqref{eq:TP_preconditioner_L20proj}, with exact inverses of the individual blocks (left), and algebraic multigrid approximations (right).}
\label{fig:cube_full_dbc}
\end{figure}

Similarly, preconditioner \eqref{eq:TP_preconditioner_L20proj} is robust in both material parameters and mesh size on the unit cube geometry, requiring slightly fewer iterations than the mixed boundary condition counterpart when applied with exact inverses (see \Cref{fig:cube_full_dbc}, left). However, we find a moderate increase in iterations with increasing $\alpha$ in the case of multigrid approximations (up to 150 iterations), while we maintain robustness with mesh refinement (see \Cref{fig:cube_full_dbc}, right). We hypothesize that this is due to the application of the SMW formula. In particular, while we control the accuracy of the AMG approximation of the operator $A^{-1}$ in \eqref{eq:smw}, the approximation of the full operator with low-rank perturbation might be significantly worse.
This suggests that a larger reduction in residual (e.g. through more effective AMG or multiple V-cycle applications) for the pressure block in preconditioner \eqref{eq:TP_preconditioner_L20proj} might be necessary to stabilize iterations in the strongly coupled regime ($\alpha \geq 1$).

\subsection{Single astrocyte test case}

\begin{figure}
    \centering
    \includegraphics[trim={3cm 0cm 3cm 6cm}, clip, width = 0.3 \textwidth]{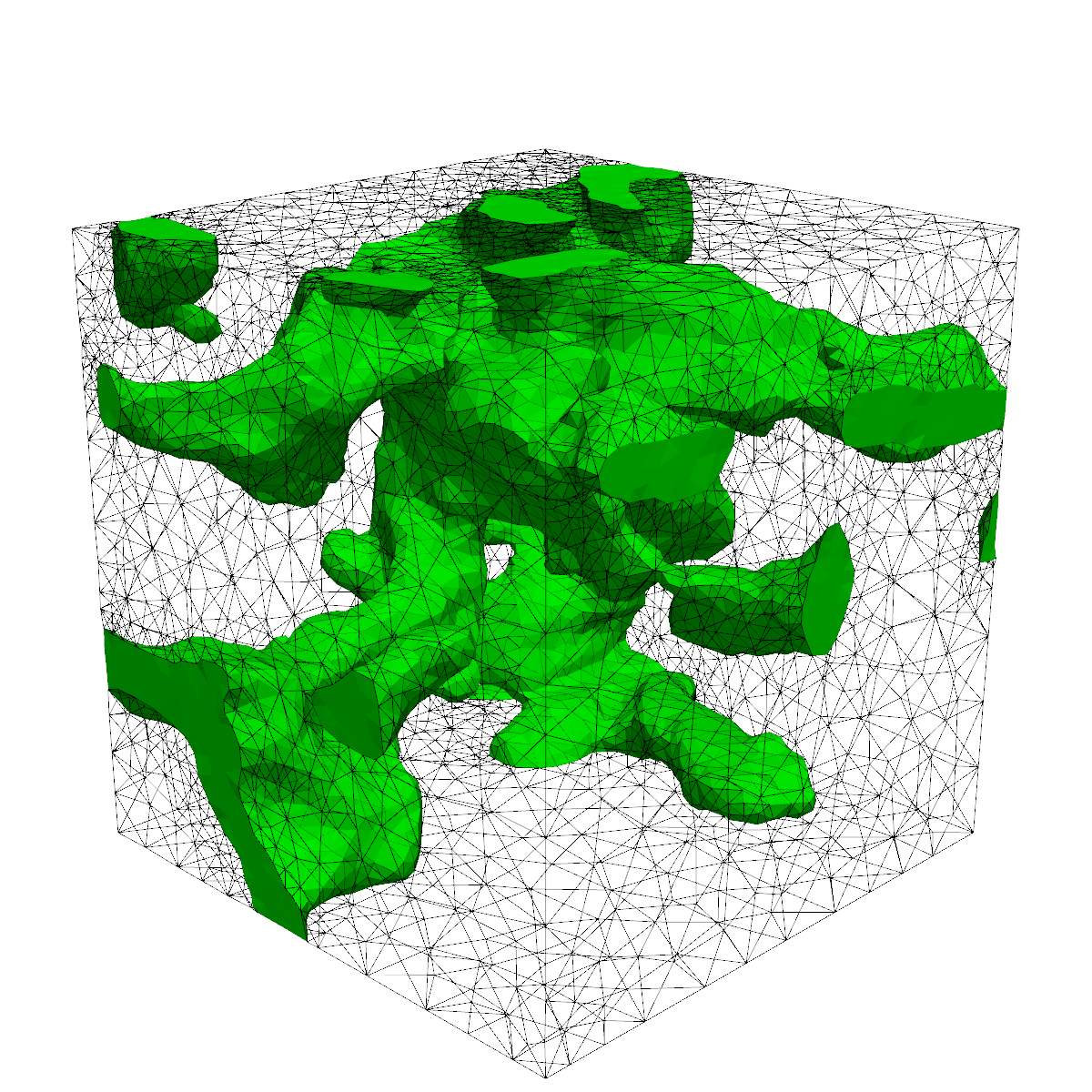}
    \includegraphics[trim={3cm 0cm 3cm 6cm}, clip, width = 0.3 \textwidth]{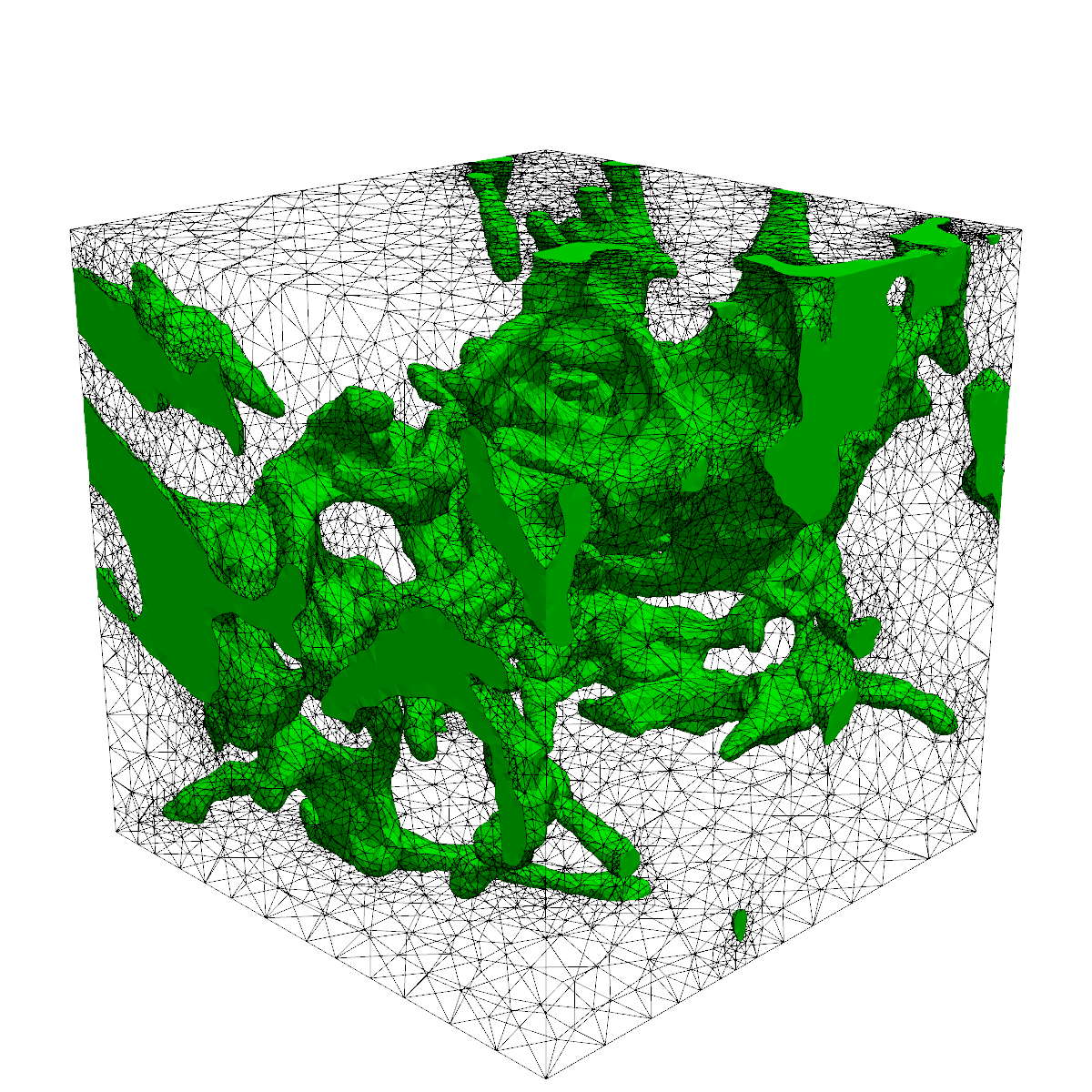}
    \includegraphics[trim={3cm 0cm 3cm 6cm}, clip, width = 0.3 \textwidth]{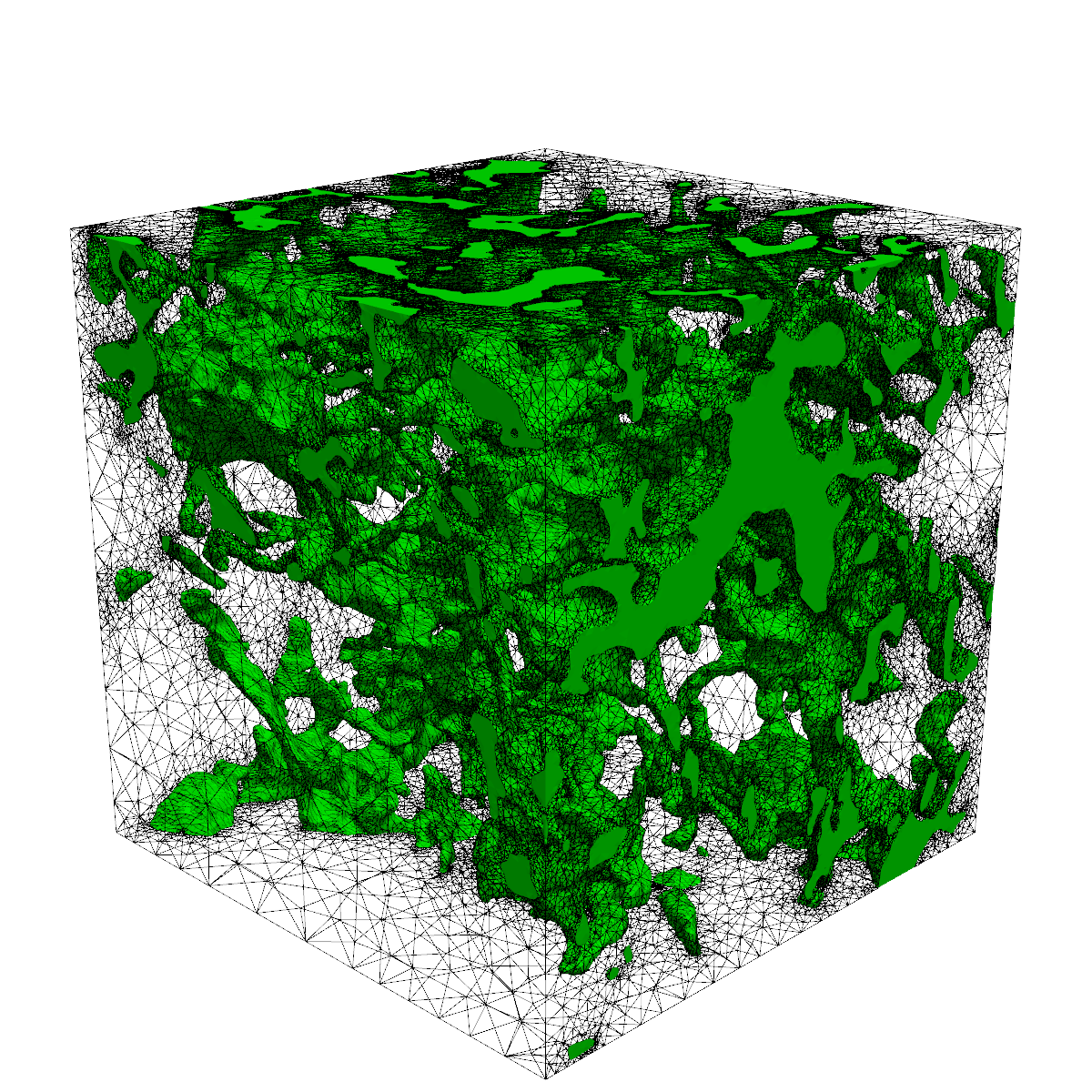}
    \caption{Sequence of meshed astrocyte geometries contained in cubes of side length $2.5$\,\textmu m, $5$\,\textmu m and $10$\,\textmu m, resulting in meshes consisting of 15\,491, 86\,538, and 606\,561 vertices, and 80\,578, 476\,286, and 3\,388\,949 cells, respectively (left to right).}
    \label{fig:astro_meshes}
\end{figure}

Aiming to evaluate the solver in a more realistic setting, we generate a sequence of astrocyte geometries based on electron microscopy images of the mouse visual cortex \cite{turner2022reconstruction}. In particular, we extract a single astrocyte contained in a cube of $2.5$\,\textmu m, $5$\,\textmu m and $10$\,\textmu m side length and generate a tetrahedral mesh of both intra- and surrounding extracellular space for each astrocyte containing cube using the meshing software fTetWild \cite{hu2020fast} (see \Cref{fig:astro_meshes}). Additionally, we refine the resulting mesh by adding an internal vertex to all cells with only boundary vertices, ensuring discrete inf-sup stability of the Taylor-Hood discretization~\cite{boffi1997three}.

We evaluate the performance of preconditioner \eqref{eq:TP_preconditioner} and \eqref{eq:TP_preconditioner_L20proj} with multigrid approximations in parallel on 64 cores and iterate until a relative tolerance of $10^{-8}$ is reached. Aiming for a realistic scenario, we scale all meshes to unit size and employ a spatially varying osmotic pressure as a right hand side
\begin{equation}
    p_{\rm osm}(x,y,z) = 1 + \sin{(2 \pi x)} \sin{(2 \pi y)} \sin{(2 \pi z)}. 
\end{equation}

\begin{figure}
    \centering
\includegraphics[trim={0.8cm 0 3cm 0},clip, height=0.48 \textwidth]{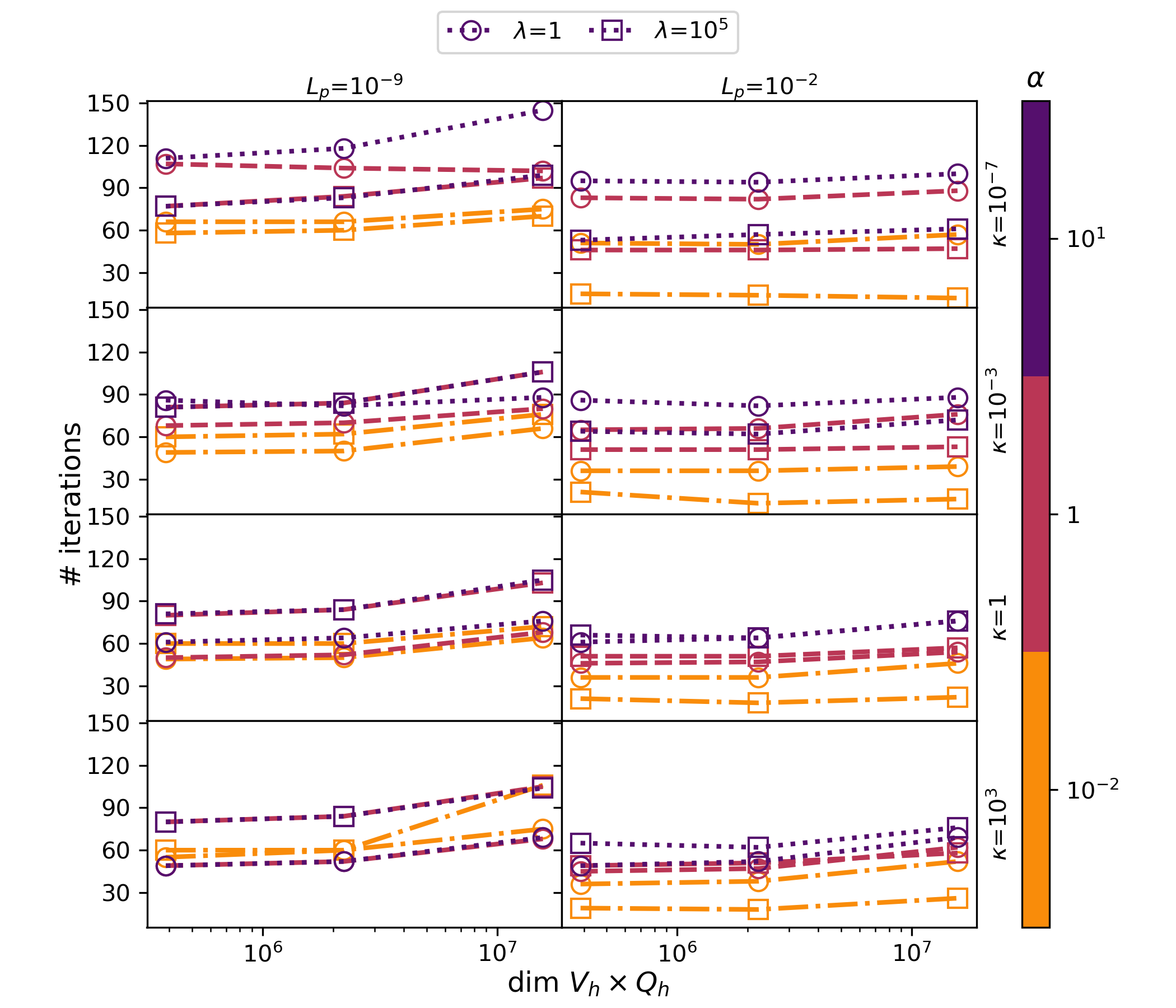}
\includegraphics[trim={0.2cm 0 0.5cm 0},clip, height=0.48 \textwidth]{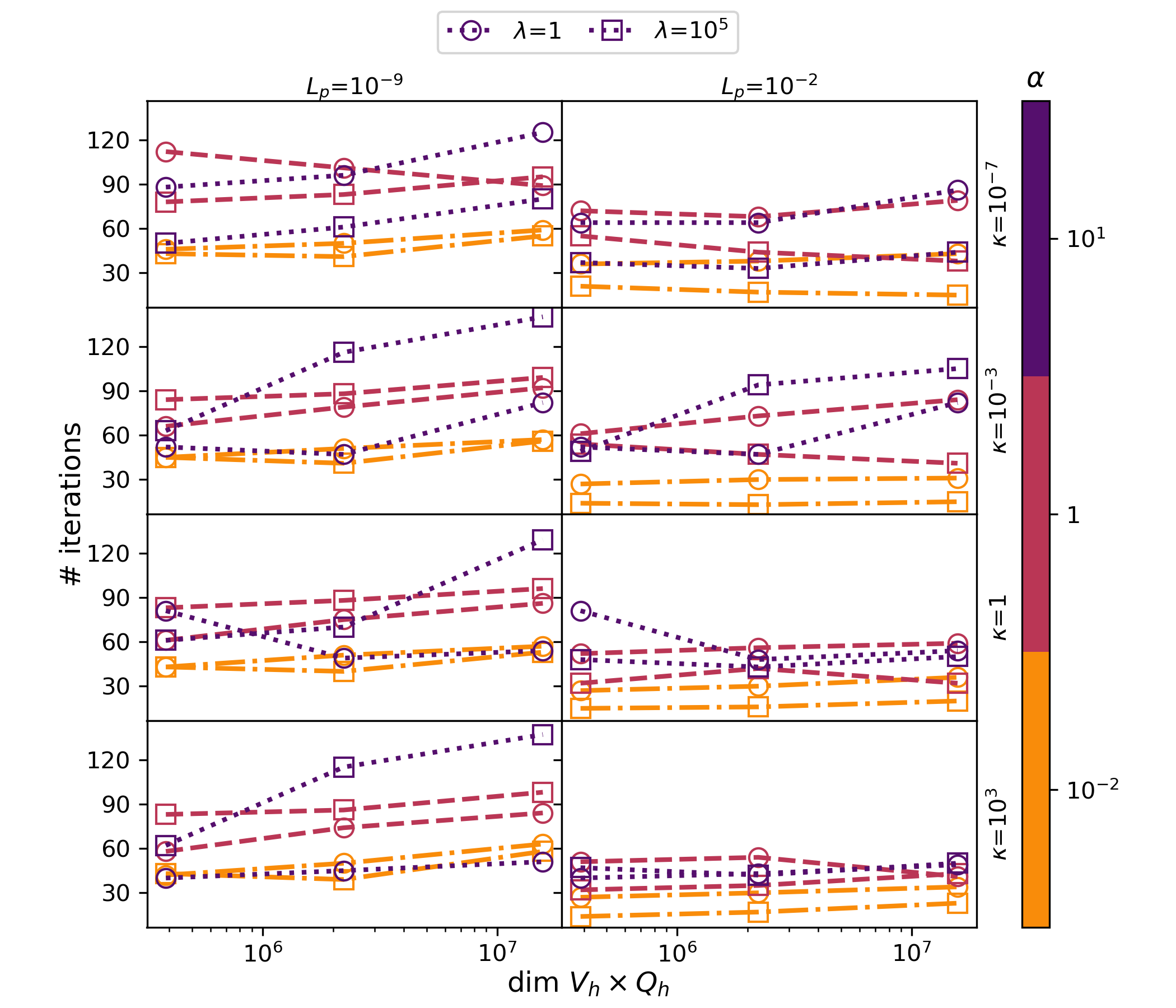}
\caption{Number of MinRes iterations for test cases with a single astrocyte with (left) mixed Dirichlet displacement boundary conditions ($\Gamma_d = \partial \Omega \cap \partial \Omega_i$) and preconditioner \eqref{eq:TP_preconditioner}, and (right) full Dirichlet boundary conditions ($\Gamma_d = \partial \Omega$) and preconditioner \eqref{eq:TP_preconditioner_L20proj}.}
\label{fig:single_astro}
\end{figure}

Starting with mixed displacement boundary conditions ($\Gamma_d = \partial \Omega \cap \partial \Omega_i$), we observe in (\Cref{fig:single_astro}, left) generally stable iteration counts, with moderate increase in some parameter regimes, mainly when the interfacial coupling is weak $L_p \ll 1$. We remark that the test setup differs from a standard mesh-refinement study; here not only the system size,
but also the geometric complexity increase substantially.

In case of full Dirichlet boundary conditions for the displacement (\Cref{fig:single_astro}, right), the application of the SMW formula in preconditioner \eqref{eq:TP_preconditioner_L20proj} can potentially amplify the multigrid approximation error, requiring higher accuracy for the approximation of the $Q$-block. In particular, we observed that insufficient accuracy of the $Q$-block AMG preconditioner yielded an indefinite operator in some cases, leading to a break-down in the MinRes solution process. For that reason, we increase the number of smoothing iterations to five, and perform two full AMG V-cycles per MinRes iteration, resulting in bounded iteration counts across all parameter variations and all system sizes.

\subsection{Cellular swelling in the mouse visual cortex}

Finally, we present an example motivated by the mechanics of cellular swelling of brain tissue. More specifically, we generated a tetrahedral mesh of a dense reconstruction of the mouse visual cortex, including the 200 largest cells in a tissue cube with 20\,\textmu m side length and an extracellular volume share of 33\,\%. The mesh consists of 4,451,225 points, with 14,336,372 and 11,259,683 extracellular and intracellular tetrahedral elements, respectively, yielding a total of 119,336,450 degrees of freedom.

\begin{figure}
    \centering
    \includegraphics[width = 1.0 \textwidth]{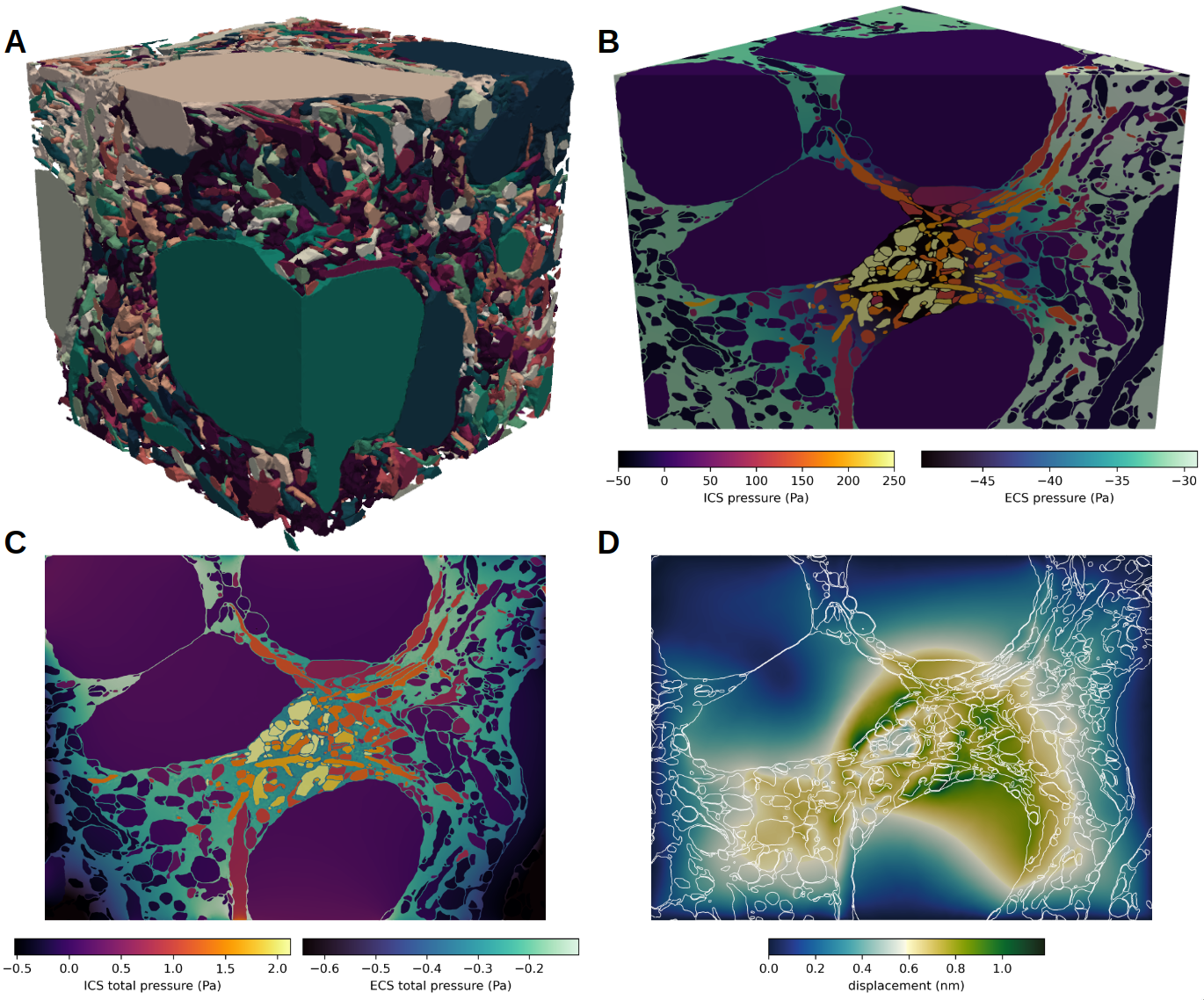}
    \label{fig:dense_swelling}
    \caption{A) Dense reconstruction of a tissue cube with 20\,\textmu m sidelength of the mouse visual cortex (200 cells, colored by cell ID);
    B) intra- and extracellular fluid pressure;
    C) intra- and extracellular total pressure;
    D) displacement magnitude in extra- and intracellular space with cell membranes (white lines).}
\end{figure}

We consider a scenario inspired by cellular swelling after intense neural activity, where an uptake of Na$^+$ and Cl$^-$ increases intracellular osmolarity, and thus induces an osmotic gradient across the cell membranes
~\cite{rungta2015cellular}. We represent the imbalance in local osmolarity ($p_{\rm osm}$) by a gaussian hill peaking at 1\,kPa and centered in the geometries midpoint, and enforce no flow and zero displacement on the outer boundary.
Choosing realistic material parameters for this setting, we assume a Young modulus of $1000\,$Pa, Poisson's ratio of $0.4$, Biot-Willis coefficient of $\alpha=1$, hydraulic conductivity of $\kappa=10^{-13}\,\text{m}^2\text{s}^{-1}\,\text{Pa}^{-1}$ and storage coefficient of $c_0=10^{-6}$\,Pa in both extra- and intracellular space, which are separated by a membrane with a permeability of $L_p=10^{-12}\,\text{m\,s}^{-1}\,\text{Pa}^{-1}$. We compute a single time step of 0.1\,s from the resting state with zero initial pressure and displacement.
Running on a single AMD EPYC 9684X node with 192 cores, the simulation consumes approximately 1.5 TB memory and successfully terminates with a relative tolerance of $10^{-8}$ after 132 iterations in 42 minutes. A visualization of the resulting pressure and displacement fields is shown in \Cref{fig:dense_swelling}, revealing a localized increase in fluid pressure in cells close to the tissue cubes midpoint (up to 250\,Pa), while the extracellular pressure drops to $-50\,$Pa. 

\section{Conclusion}

We have presented efficient and robust solvers for a cell-by-cell poroelasticity model motivated by cellular mechanics in brain tissue. Leveraging the framework of norm-equivalent preconditioning \cite{mardal_preconditioning_2011} and fitted norms \cite{hong2023new}, we have established the preconditioners robustness with respect to variations in material parameters, and demonstrated its performance across a large range of parameters relevant to realistic scenarios through numerical experiments. The method's scalability to large and complex geometries is achieved through efficient AMG realizations of these preconditioners, ensuring computationally inexpensive yet spectrally equivalent operator approximations. Importantly, we covered a broad set of boundary conditions, including full Dirichlet boundary conditions for the displacement, where we efficiently handled the arising dense projection operator through the SMW formula.
Future research will focus on applying this methodology to investigate cellular swelling and volume regulation, leveraging the model's capability to simulate these complex physiological processes in detailed cellular geometries.

%\section*{Acknowledgments}

\bibliographystyle{siamplain}
\bibliography{references}
\end{document}

%% file: ex_shared.tex
\usepackage{lipsum}
\usepackage{amsfonts}
\usepackage{graphicx}
\usepackage{epstopdf}
\usepackage{algorithmic}
\usepackage{enumitem}
\usepackage{nicematrix, tikz}
\usepackage{float}
\usepackage{a4wide}
\usepackage{makecell}
\usepackage{diagbox}
\usepackage{amssymb, stmaryrd}
\usepackage{booktabs}
\usepackage{siunitx}
\usepackage{cleveref}

\ifpdf
  \DeclareGraphicsExtensions{.eps,.pdf,.png,.jpg}
\else
  \DeclareGraphicsExtensions{.eps}
\fi

\newsiamremark{remark}{Remark}
\newsiamremark{hypothesis}{Hypothesis}
\newsiamremark{example}{Example}
\crefname{hypothesis}{Hypothesis}{Hypotheses}
\newsiamthm{claim}{Claim}
\newsiamremark{fact}{Fact}
\crefname{fact}{Fact}{Facts}

\headers{Robust solver for a cell-by-cell poroelasticity problem}{M. Causemann and M. Kuchta}

\title{Parameter-robust preconditioners for a cell-by-cell poroelasticity model with interface coupling\thanks{Submitted to the editors 26.05.2025.
    \funding{This project received from funding from the Research Council of Norway via grants \#324239 (EMIx),
      \#303362 (DataSim) and the European Research Council via grant \#101141807 (aCleanBrain).}}}

\author{Marius Causemann\thanks{Simula Research Laboratory, Oslo, Norway 
  (\email{mariusca@simula.no})}
\and Miroslav Kuchta\thanks{Simula Research Laboratory, Oslo, Norway 
  (\email{miroslav@simula.no})}}

\usepackage{amsopn}

\newcommand{\eps}{\ensuremath{\mathbf{\epsilon}}}
\renewcommand{\d}{\ensuremath{\mathbf{d}}}

\renewcommand{\dh}{\ensuremath{\mathbf{d}_h}}
\newcommand{\I}{\ensuremath{\mathbf{I}}}

\newcommand{\pfh}{p_{F,h}}
\newcommand{\pth}{p_{T,h}}

\renewcommand{\v}{\ensuremath{\mathbf{v}}}

\renewcommand{\u}{\ensuremath{\mathbf{u}}}

\renewcommand{\div}{\text{ div }}
\newcommand{\f}{\ensuremath{\mathbf{f}}}

\newcommand{\n}{\ensuremath{\mathbf{n}}}

\newcommand{\jump}[1]{\llbracket #1 \rrbracket}

\NiceMatrixOptions
  {
    custom-line = 
     {
       letter = I ,
       command = dashedline ,
       tikz = densely dashed ,
     }
  }

  \crefname{subsection}{section}{sections}
  \Crefname{subsection}{Section}{Sections}